\def\versiondate{23 Jan.\ 2003}
\input math.macros
\input lrlEPSfig.macros

\checkdefinedreferencetrue
\continuousfigurenumberingtrue
\theoremcountingtrue
\sectionnumberstrue
\forwardreferencetrue
\citationgenerationtrue
\nobracketcittrue
\tocgenerationtrue
\hyperstrue
\initialeqmacro

\bibsty{myapalike}

\def\T{{\Bbb T}}
\def\Td{{\T^d}}
\def\dbar{{\overline{d}}}
\def\Pfp{{(\P^f)^+}}
\def\Pfm{{(\P^f)^-}}
\def\calA{{\cal A}}
\def\calT{{\cal T}}
\def\F{{\cal F}}

\def\C{{\Bbb C}}
\def\hE{\ell^2(E)}
\def\Hs{{\ell^2(E)}}
\def\hD{\ell^2({\Z}^d)}

\def\ip#1{(\changecomma #1)}
\def\bigip#1{\bigl(\bigchangecomma #1\bigr)}

\def\changecomma#1,{#1,\,}
\def\bigchangecomma#1,{#1,\;}
\def\leftchangecomma#1,{#1,\ }

\def\bfz{{\bf 0}}
\def\bfo{{\bf 1}}
\def\GM{{\ss GM}}  
\def\mult#1{M_{\times #1}} 
\def\strle{\preccurlyeq_{\rm s}}   
\def\fav{\widehat f(0)}  
\def\past{{\ss Past}}
\def\cat{{\rm G}}  
\def\others{{\cal Z}}  
\def\HM{{\ss HM}}  
\def\fullle{\preccurlyeq_{\rm f}}   
\def\Zd{{\Z^d}}
\def\ee{{\bf e}}  
\def\frac#1#2{{#1 \over #2}}
\def\nonfrac#1#2{{#1 / #2}}
\def\hl{{\ss HL}}  
\def\supp{{\rm supp}\,}
\def\Cov{{\rm Cov}}
\def\ord{\Pi}  

\def\BLPSusf{Benjamini, Lyons, Peres, and Schramm (2000)%
\def\BLPSusf{BLPS (2000)}}

\def\firstheader{\eightpoint\ss\underbar{\raise2pt\line 
    {
    To appear in {\it Duke Math. J.}
    \hfil Version of \versiondate}}}

\beginniceheadline

\vglue20pt

\title{Stationary Determinantal Processes:}
\title{Phase Multiplicity, Bernoullicity,}
\title{Entropy, and Domination}

\author{Russell Lyons and Jeffrey E. Steif}

\abstract{%
We study a class of stationary processes indexed by $\Z^d$ that
are defined via minors of $d$-dimensional (multilevel) Toeplitz matrices.
We obtain necessary and sufficient conditions for
phase multiplicity (the existence of a phase transition)
analogous to that which occurs in 
statistical mechanics.
Phase uniqueness is equivalent to the presence of a strong
$K$ property, a particular strengthening of the usual $K$ 
(Kolmogorov) property. 
We show that all of these processes are
Bernoulli shifts (isomorphic to i.i.d.\ processes in the sense of ergodic
theory).
We obtain estimates of their entropies and
we relate these processes via stochastic domination to product measures.
}

\bottomII{Primary 
82B26, 
28D05, 
60G10. 
Secondary
82B20, 
37A05, 
37A25, 
37A60, 
60G25, 
60G60, 
60B15. 
}
{Random field, Kolmogorov mixing, Toeplitz determinants, negative
association, fermionic lattice model, stochastic domination, entropy,
Szeg\H{o} infimum, Fourier coefficients, insertion tolerance, deletion
tolerance, geometric mean.}
{Research partially supported by 
NSF grants DMS-9802663, DMS-0103897 (Lyons), and
DMS-0103841 (Steif), the Swedish Natural Science Research Council (Steif)
and the Erwin Schr\"odinger Institute for Mathematical
Physics, Vienna (Lyons and Steif).}

\articletoc

\bsection {Introduction}{s.intro}

Determinantal probability measures and point processes arise in numerous
settings, such as
mathematical physics (where they are called fermionic point processes),
random matrix theory, representation theory, and certain other areas
of probability theory. 
See \ref b.Sosh/ for a survey and \ref b.L:det/ for additional developments
in the discrete case.
We present here a detailed analysis of the discrete
stationary case.
After this paper was first written, we learned of independent but slightly
prior work of \ref b.ShiTak:II/, announced in \ref b.ShiTak:ann/.
We discuss the (small) overlap between their work and ours in the appropriate
places below.
See also \ref b.ShiTak:I/ and \ref b.ShiYoo:glauber/ for related
contemporaneous work.

Stationary determinantal processes are interesting from several viewpoints.
First, they have interesting relations with the 
theory of Toeplitz determinants.
As in that theory, the {\bf geometric mean} of a nonnegative function $f$,
defined as 
$$
\GM(f) := \exp \int \log f
\,,
$$
will play an important role in some of our results.
(In fact, the arithmetic mean and the harmonic mean of $f$ will also
characterize certain properties of our processes.)
Second, such processes arise in certain combinatorial models, such as uniform
spanning trees and dimer models.
Third, these systems have a rich infinite-dimensional parameter space,
consisting, in the case of a $\Z^d$ action, of all measurable functions $f$
from the $d$-dimensional torus to $[0, 1]$.
We illustrate the variety of possible behaviors with examples throughout
the paper.
Fourth, they have the unusual property of negative association.
Though unusual, negative association occurs in various places and a fair
amount is known about it
(see \ref b.Joag/, \ref b.Pemantle:negass/, 
\ref b.Newman:asympt/, \ref b.ShaoSu:LIL/, \ref b.Shao:comparison/,
\ref b.ZhangWen:weak/, \ref b.Zhang:Strassen/, and the references therein,
for example). We offer a whole new class of examples of negatively
associated stationary processes.
As such, determinantal processes provide an easy way to construct examples of
many kinds of behavior that might otherwise be difficult to construct, such as
negatively associated (stationary) processes with slow decay of correlations,
or with the even sites independent of the odd sites (in one dimension, say),
or with the property of being finitely dependent.
Fifth, all our processes are Bernoulli shifts, i.e., isomorphic to i.i.d.\
processes. This may be surprising in view of the fact that only measurability,
rather than smoothness, of the parameter $f$ is required.
Sixth, in one dimension,
some of the processes are strong $K$, while others are not.
Namely, strong $K$ is equivalent to $f(\bfo-f)$ having a positive
geometric mean.
Similarly, we characterize exactly, in all dimensions, which $f$ give strong
{\it full\/} $K$ systems. It turns out that not only the rate at which $f$
approaches 0 or 1 matters, but also where. For example, in two dimensions,
if $f$ is real analytic, then the system is strong full $K$ iff the
(possibly empty) sets $f^{-1}(0)$ and $f^{-1}(1)$ belong to nontrivial
algebraic varieties. The strong full $K$ property is 
analogous to phase uniqueness
in statistical physics, as
we explain in \ref s.sk/.

We now state our results somewhat more precisely and present several
examples.
Let $f : \Td \to [0, 1]$ be a Lebesgue-measurable function on the
$d$-dimensional torus $\Td := \R^d/\Z^d$.
Define a $\Z^d$-invariant probability measure $\P^f$ on the Borel sets of $ \{
0, 1 \}^{\Z^d}$ by defining the probability of the cylinder sets 
\begineqalno
\P^f[\eta(e_1)=1,\dots, \eta(e_k)=1] 
&:= 
\P^f[\{ \eta \in \{ 0, 1 \}^{\Z^d} \st \eta(e_1)=1,\dots, \eta(e_k)=1\}]
\cr&:= 
\det [\widehat{f}(e_j-e_i)]_{1\le i, j \le k}
\cr
\endeqalno
for all $e_1,\ldots,e_k\in \Z^d$, where $\widehat f$ denotes the Fourier
coefficients of $f$.
We shall prove in \ref s.back/ that this does indeed define a probability
measure.
Note that when $d=1$ and when $e_1,\ldots,e_k$
are chosen to be $k$ consecutive integers, the
right-hand side above is the usual $k \times k$ Toeplitz determinant of $f$
denoted $D_{k-1}(f)$.
In particular, we have $\P^f\big[\eta(e) = 1\big] = \widehat f(\bfz)
= \int_{\Td} f$ for every $e \in \Z^d$.

\procl x.Bernoulli
As a simple example, if $f \equiv p$, then $\P^f$ is i.i.d.\ Bernoulli($p$)
measure.
\endprocl

\procl x.ustH
For a more interesting example, consider
$$
f(x, y) := {\sin^2 \pi x \over \sin^2 \pi x + \sin^2 \pi y}
\,.
\label e.ustH
$$
A portion of a sample of a configuration from $\P^f$ is shown in \ref
f.window-ustH/, where a square with upper left corner $(i, j)$ is colored
black iff $\eta(i, j) = 1$.
These correspond to the horizontal edges of a uniform spanning tree in the
square lattice. That is, if $T$ is a spanning tree of $\Z^2$, then $\eta(i,
j)$ is the indicator that the edge from $(i, j)$ to $(i+1, j)$ belongs to
$T$. The portion of the spanning tree from which \ref f.window-ustH/ was
constructed is shown in \ref f.window-ust/.
When $T$ is chosen ``uniformly" (see \ref b.Pemantle:ust/, \ref
b.Lyons:bird/, or \ref b.BLPSust/ for definitions and information on
this), then $\eta$ has the law $\P^f$. This follows from the Transfer
Current Theorem of \ref b.BurPem/ and the representation of the Green
function as an integral; see \ref b.Lyons:book/ for more details.
Similarly, the edges of the uniform spanning forest (it is a tree only for $d
\le 4$, as shown by \ref b.Pemantle:ust/) parallel to the $x_1$-axis in $d$
dimensions correspond to the
function 
$$
f(x_1, x_2, \ldots, x_d) := {\sin^2 \pi x_1 \over \sum_{j=1}^d \sin^2 \pi x_j}
\,.
\label e.ustH-d
$$
We remark that
the uniform spanning tree is the so-called random-cluster model when 
one takes the limit $q \downarrow 0$,
then $p \downarrow 0$, and finally the thermodynamic limit, the latter
shown to exist by \ref b.Pemantle:ust/.
\endprocl

\twoefiginslabel window-ustH {A sample from $\P^f$ of \ref e.ustH/.} x3 
window-ust {A uniform spanning tree.} y3

\procl x.ustHx
Let 
$$
g(x) :=
{\sin \pi x \over \sqrt{1+\sin^2 \pi x}}
\,.
$$
Then the edges of the uniform spanning tree in the plane that
lie on the $x$-axis have the law $\P^g$, as shown in \ref x.ustH_x-dom/ below.
\endprocl

\procl x.USTzz
Let 
$$
f(x) := {1 \over 2} + {|\sin 2\pi x| - 1 \over 2 \cos 2\pi x}
\,.
$$
An elementary calculation shows that 
$$
\widehat f(k) = \cases{1/2  &if $k=0$,\cr
                       0    &if $k \ne 0$ is even,\cr
\displaystyle
(-1)^{(k-1)/2}\left(-{1 \over 2} + {2 \over \pi} 
      \sum_{j=0}^{(k-1)/2} {(-1)^j \over 2j+1}\right) - {1 \over \pi k}
                            &if $k$ is odd.\cr}
$$
We shall show in \ref x.USTzz-dom/ that the measure $\P^f$ arises as follows.
Given a spanning tree $T$ of the square lattice, let $\eta(n)$ be the
indicator that $e_n \in T$, where $e_n$ is the edge 
$$
e_n := \cases{[(n/2, n/2), (n/2 + 1, n/2)]  &if $n$ is even,\cr
              [((n+1)/2, (n-1)/2), ((n+1)/2, (n+1)/2)]  &if $n$ is odd.\cr}
$$
The collection of edges $\{e_n \st n \in \Z \}$ is a zig-zag path in the
plane. The law of $\eta$ is $\P^f$ when $T$ is chosen as a uniform spanning
tree.
(Although it is not hard to see that the law of $\eta$ is $\Z$-invariant,
by using planar duality, and although the law of $\eta$ must be a
determinantal probability measure because the law of $T$ is,
it is not apparent {\it a priori\/} that the law of $\eta$ has the
form $\P^f$ for some $f$.)
\endprocl

\procl x.lozenge
Fix a horizontal edge of the hexagonal lattice (also known as the
honeycomb lattice) and index all its vertical translates by $\Z$.
If one considers the standard measure of maximal entropy on perfect
matchings of the hexagonal lattice, also called the dimer model
and equivalent to lozenge tilings of the plane, and
looks only at the edges indexed as above by $\Z$, then the law is $\P^f$
for $f := \I{[1/3, 2/3]}$, as shown by \ref b.Kenyon:local/.
\endprocl

\procl x.domino It is interesting that the function
$f := \I{[0, 1/2]}$ for $d = 1$ also arises from a combinatorial model.
In this case,
$$
\widehat f(n) = \cases{1/2 &if $n=0$,\cr
                  0 &if $n \ne 0$ is even,\cr
                  1/(\pi i n) &if $n$ is odd.\cr}
$$
The measure $\P^f$ is the zig-zag process of \ref b.Joh:noninter/
derived from uniform domino tilings in the plane.
For the definition of ``uniform" in this case, see \ref b.BurPem/.
A picture of a portion of such a tiling is shown in \ref f.window-tiling/.
Consider the squares on a diagonal from upper left to lower right.
The domino covering any such square also covers a second square.
If this second square is above the diagonal, then we color the original
square black, as shown in \ref f.window-diagprocess/.
\ref b.Joh:noninter/ showed that the law of this process is $\P^f$ when the
diagonal squares are indexed by $\Z$ in the natural way.
More generally, the processes $\P^f$ for $f$ the indicator of any arc of $\T$
are used by \ref b.BOO:plancherel/, Theorem 3, to describe the typical shape
of Young diagrams.
\endprocl

\twoefiginslabel window-tiling {A uniform domino tiling.} x3 
window-diagprocess {A sample from $\P^f$ of \ref x.domino/.} y3

\procl x.renewal Let $0 < a < 1$ and $d = 1$.
If $f(x) := (1 - a)^2/|e^{2\pi i x} - a|^2$, then $\P^f$ is a renewal
process (\ref B.Sosh/).
The number of 0s between successive 1s has the same distribution as
the number of tails until 2 heads appear for a coin that has probability $a$
of coming up tails.
More explicitly, for $n\ge 1$,
$$
\P^f[\eta(1) = \cdots = \eta(n-1) = 0, \eta(n) = 1 \mid \eta(0) = 1]
=
n (1-a)^2 a^{n-1}
\,.
\label e.renewal-distribution
$$
Since 
$$
f(x) = {1-a \over 1+a} \left( {a e^{2\pi i x} \over 1 - a e^{2\pi i x}} + {1
\over 1 - a e^{-2\pi i x}} \right)
\,,
$$
expansion in a geometric series shows that
$$
\widehat f(k) = {1-a \over 1+a} a^{|k|}
\,.
$$
We prove that $\P^f$ is indeed this explicit renewal process after we prove
\ref p.regenerate/, in which we extend this example to other regenerative
processes.
\endprocl

\procl x.perconf
If $0 < p < 1$ and $f : \Td \to [0, 1]$ is measurable, then a fair sample of
$\P^{pf}$ can be obtained from a fair sample of $\P^f$ simply by independently
changing each 1 to a 0 with probability $1-p$.
\endprocl

For general systems, note that the covariance of
$\eta(\bfz)$ and $\eta(k)$ for $k \in \Z^d$ is $-|\widehat f(k)|^2$. This is
summable in $k$ since $f \in L^2(\Td)$, but that is essentially the most one
can say for its rate of decay.
That is, given any $\Seq{a_k} \in \ell^1(\Z^d)$, there is some
(even continuous) $f : \Td \to [0, 1]$ and some constant $c > 0$ such that
$|\widehat f(k)|^2 \ge c |a_k|$ for all $k \in \Z^d$, as shown by \ref
b.LKK/.
Observe also that as is the case for Gaussian processes, the processes
studied here have the property that if the random variables are uncorrelated, then
they are mutually independent.

It is shown (in a much more general context) by \ref b.L:det/ that these
measures, as well as these measures conditioned on the values of $\eta$
restricted to any finite subset of $\Z^d$,
have the following negative association property:
If $A$ and $B$ are increasing events that are measurable with respect to the
values of $\eta$ on disjoint subsets of $\Z^d$, then $A$ and $B$ are
negatively correlated.

In the order presented in this paper, our principal results are the following.
\beginbullets

For all $f$, the process $\P^f$ is a Bernoulli shift. 
This was shown by \ref b.ShiTak:II/ for those $f$ such that $\sum_{n \ge 1}
n |\widehat f(n)|^2 <\infty$ by showing that those $\P^f$ are weak Bernoulli.

For all $f$, the process $\P^f$ stochastically dominates product measure
$\P^{\GM(f)}$ and is stochastically dominated by product measure
$\P^{1-\GM(\bfo-f)}$, and these bounds are optimal.
This is rather unexpected for the process of \ref x.ustH/ related to the
uniform spanning tree; explicit calculations are given below in \ref
x.ustH-dom/ for this particular process.
We give similar optimal bounds for full domination (uniform insertion and
deletion tolerance) in terms of harmonic means.

We present methods to estimate the entropy of $\P^f$.
For example, we show that for the function $g$ in \ref x.ustHx/,
the entropy of the system $\P^g$ lies in the interval $[0.69005, 0.69013]$;
see \ref x.ustH-x-ent/.

The process $\P^f$ is strong full $K$ iff there is a nonzero
trigonometric polynomial $T$ such that 
${|T|^2 \over f (\bfo-f)} \in L^1(\Td)$.
This is equivalent to phase uniqueness in the sense that
no conditioning at infinity can change the measure.

In one dimension, $\P^f$ is strong $K$ iff $f(\bfo-f)$ has a positive
geometric mean.
This is equivalent to phase uniqueness when conditioning on one
side only. Higher-dimensional versions of this will also be obtained.
\endbullets

We shall give full definitions as they become needed. 
Some general background on determinantal probability measures is presented
in \ref s.back/, where we also exhibit a key representation of certain
conditional probabilities as Szeg\H{o} infima.
The property of being a Bernoulli shift 
is proved in \ref s.bern/, while the auxiliary result
that all $\P^f$ have full support (except in two degenerate cases) is shown in
\ref s.fs/.
Properties concerning stochastic domination are proved in \ref
s.dom/. These are used to estimate entropy in \ref s.ent/.
More sophisticated methods
of estimating entropy are also developed and illustrated in \ref s.ent/.
The main result about phase multiplicity is proved in \ref s.sk/, while the
one-sided case for $d=1$ and its higher-dimensional generalizations are
treated in \ref s.1side/.
Finally, we end with some open questions in \ref s.open/.

The above definitions can be generalized to any countable abelian group with
the discrete topology. 
For example, if we use $\Z \times \Z_2$, 
we obtain nontrivial joinings of the above systems with themselves. That
is, suppose $f : \T \to [0, 1]$ and $h : \T \to [0, 1]$ are such that $f \pm
h : \T \to [0, 1]$. Then the function $f_h : \T \times \{ -1, 1 \}
\to [0, 1]$ defined by $(x, \epsilon) \mapsto f(x) + \epsilon h(x)$
gives a system that, when restricted to each copy of $\Z$, is just $\P^f$, but
has correlations between the two copies that are given by $h$
(the case $h = \bfz$ gives the independent joining). 
We can obtain slow decay of correlations between the two copies and
negative associations in the joining itself;
this is perhaps something that is not easy to construct directly.

\bsection {Background}{s.back}

We first quickly review the probability measures studied in 
\ref b.L:det/;
see that paper for complete details.

Let $E$ be a finite or countable set and consider the complex
Hilbert space $\hE$. 
Given any closed subspace $H \subseteq \ell^2(E)$,
let $P_H$ denote the orthogonal projection onto $H$.
There is a unique probability measure $\P^H$ on $2^E := \{ 0, 1 \}^E$ defined
by 
$$
\P^H[\eta(e_1) = 1, \ldots, \eta(e_k) = 1]
=
\det [(P_H e_i, e_j)]_{1\le i,j\le k}
\label e.(3)
$$
for all $k \ge 1$ and any set of distinct $e_1,\ldots,e_k\in E$; see, e.g.,
\ref b.L:det/ or \ref b.DVJ/, Exercises 5.4.7--5.4.8.
On the right-hand side, we are identifying each $e \in E$ with the element
of $\hE$ that is 1 in coordinate $e$ and 0 elsewhere. 
In case $H$ is finite-dimensional, then $\P^H$ is concentrated on subsets
of $E$ of cardinality equal to the dimension of $H$.

More generally,
let $Q$ be a positive contraction, meaning that $Q$ is a self-adjoint
operator on $\Hs$ such that for all $u \in \Hs$, we have $0 \le (Q u, u) \le
(u, u)$.
There is a unique probability
measure $\P^Q$ such that 
$$
\P^Q[\eta(e_1) = 1, \ldots, \eta(e_k) = 1]
= \det [(Q {e_i}, {e_j})]_{i, j \le
k}
\label e.Qcorrel
$$
for all $k \ge 1$ and distinct $e_1, \ldots, e_k \in E$.
When $Q$ is the orthogonal projection onto a closed
subspace $H$, then $\P^Q = \P^H$.
In fact, properties of $\P^Q$ can be deduced from the special case of
orthogonal projections. Since this will be useful for our analysis, we
review this reduction procedure.

Note first that uniqueness follows from the fact that \ref e.Qcorrel/
determines all finite-dimen\-sion\-al marginals via the inclusion-exclusion
theorem. Indeed, we have the following formula for any disjoint pair of
finite sets $A, B \subseteq E$ (see, e.g., \ref b.L:det/):
$$
\P^Q[\eta \restrict A \equiv 1,\, \eta \restrict B \equiv 0]
=
\det \left[\big(\I B(e)e + (-1)^{\I B(e)} Q e, e'\big)\right]_{e, e' \in A
\cup B}
\,.
\label e.incexc
$$
To show existence, let $P_H$ be any orthogonal projection that is
a dilation of $Q$, i.e., $H$ is
a closed subspace of $\ell^2(E')$ for some $E' \supseteq E$ and for all $u
\in \Hs$, we have $Qu = P_{\Hs} P_H u$, where we regard $\ell^2(E') = \Hs
\oplus \ell^2(E' \setminus E)$. (In this case, $Q$ is also called the
compression of $P_H$ to $\Hs$.)
The existence of a dilation is standard and is easily
constructed; see, e.g., \ref b.L:det/.
Having chosen a dilation, we simply define $\P^Q$ as the law of $\eta$
restricted to $E$ when $\eta$ has the law $\P^H$.
Then \ref e.Qcorrel/ is a special case of \ref e.(3)/.

A probability measure $\P$ on $2^E$ is said to have {\bf negative
associations} if for all pairs
$A$ and $B$ of increasing events that are measurable with respect to the
values of $\eta$ on disjoint subsets of $E$, we have that $A$ and $B$ are
negatively correlated with respect to $\P$.
The following {\bf conditional negative association} (CNA)
property is proved in \ref b.L:det/ and a consequence of this
(see \ref p.JNRD/) will be used frequently here.

\procl t.negass
If $Q$ is any positive contraction on $\ell^2(E)$, $A$ is a
finite subset of $E$, and $\eta_0 \in 2^A$,
then $\P^Q[\;\cdot \mid \eta \restrict A = \eta_0]$ has
negative associations.
\endprocl

We now assume that $E=\Z^d$. Then the group structure of $\Z^d$ allows
$\Z^d$ to act naturally on $\ell^2(\Z^d)$ and on $2^{\Z^d}$.
The proof of the following lemma is straightforward and 
therefore skipped.

\procl l.pm
If $Q$ is a $\Z^d$-invariant positive contraction on
$\ell^2(\Z^d)$, then $\P^Q$ is also $\Z^d$-invariant.
\endprocl

As is well known, there exists a complex Hilbert-space isomorphism
between $L^2(\T^d, \lambda_d)$ and $\ell^2(\Z^d)$ where
$\T^d$ is the $d$-dimensional torus $\R^d/\Z^d$ and
$\lambda_d$ is unit Lebesgue measure on $\T^d$. 
This isomorphism is given by the Fourier transform
$f\mapsto \widehat f$, where for $f \in L^2(\T^d, \lambda_d)$, we have
$\widehat{f}(k)=\int_{\T^d} f(x) e^{-2\pi i k\cdot x} \, d\lambda_d(x)$ 
for $k\in\Z^d$.
If $\ee_k$ denotes the function $x \mapsto e^{2\pi i k \cdot x}$, then the
isomorphism takes the set $ \{ \ee_k \st k \in \Zd \}$ to the standard basis
for $\ell^2(\Zd)$.
From now on, we shall abbreviate $L^2(\T^d, \lambda_d)$ by $L^2(\T^d)$.

The following is well known.

\procl t.rudin 
\beginitems
\itemrm{(i)} Let $A\subseteq\Td$ be measurable and consider the operator 
$T_A:L^2(\T^d)\to L^2(\T^d)$ given by
$$
T_A(g)=g \I A\,.
$$
Then these projections (as $A$ varies over the measurable subsets of $\Td$) 
correspond (via the Fourier
isomorphism) to the $\Z^d$-invariant projections on $\ell^2(\Z^d)$.
\itemrm{(ii)} More generally, let $f:\Td\to [0,1]$ be measurable and
consider the operator $M_f:L^2(\T^d)\to L^2(\T^d)$ given by
$$
M_f(g)=fg\,.
$$
Then these positive contractions (as $f$ varies
over the measurable functions from $\T^d$ to $[0,1]$) correspond (via the
Fourier
isomorphism) to the $\Z^d$-invariant positive contractions on $\ell^2(\Z^d)$.
More specifically, $M_f$ corresponds to convolution with $\widehat{f}$.
\enditems
\endprocl

As in the above theorem, an $f:\Td\to [0,1]$ yields a
$\Z^d$-invariant positive contraction $Q_f$ on
$\ell^2(\Z^d)$, which in turn yields a translation-invariant
probability measure $\P^{Q_f}$ on $2^{\Z^d}$ that we denote more simply
by $\P^f$.

\procl l.fourier
Given $f:\Td\to [0,1]$ measurable and $e_1,\ldots,e_k\in \Z^d$,
$$
\P^f[\eta(e_1)=1,\dots, \eta(e_k)=1] = 
\det [\widehat{f}(e_j-e_i)]_{1\le i, j \le k}
\,.
$$
\endprocl

\proof 
By definition, the left-hand side is
$\det [(Q_f e_i,e_j)]_{1\le i, j \le k}$. 
By \ref t.rudin/(ii), 
$$
\ip{Q_f e_i,e_j}= \bigip{M_f \ee_{e_i}, \ee_{e_j}}
= \widehat{f}(e_j-e_i)\,.
\Qed
$$

\procl r.toeplitz
\ref l.fourier/ says that for $d=1$,
the probability of having 1s on some finite collection of elements of
$\Z$ is a particular minor of the Toeplitz matrix associated to $f$.
\endprocl

Equation \ref e.incexc/ shows a symmetry of $\P^f$ and
$\P^{\bfo-f}$, namely, if $\eta$ has the distribution $\P^f$, then $\bfo - \eta$
has the distribution $\P^{\bfo-f}$.
\ref b.ShiTak:II/ prove the existence of $\P^f$ by a different method and
also note this symmetry.

Although the last lemma gives us a formula for $\P^f$ directly in terms of
$f$ without reference to any projections, it is still useful 
to know a specific projection of which $M_f$ is a compression. 
Let 
$f:\Td\to [0,1]$ be measurable. 
Identifying $\T^{d+1}$ with $\Td \times \CO{0, 1}$,
we let $A_f\subseteq \T^{d+1}$ be the set
$\{(x,y)\in \T^{d+1} \st y\le f(x)\}$. 
Consider the projection $T_{A_f}$ of $L^2(\T^{d+1})$ given 
in \ref t.rudin/(i).
We view $L^2(\T^d)$ as a subspace of $L^2(\T^{d+1})$ by identifying
$g\in L^2(\T^d)$ with $g\otimes \bfo \in L^2(\T^{d+1})$, where
$(g\otimes \bfo) (x,y):= g(x)$ for $x\in \Td$, $y\in \T$.
The orthogonal projection $P$ of $L^2(\T^{d+1})$ onto $L^2(\T^d)$ is then
given by $g \mapsto (x \mapsto \int_\T g(x,y) \, dy)$. A very simple
calculation, 
left to the reader, shows that $M_f$ is a compression of $T_{A_f}$;
i.e., 
$$
M_f=P T_{A_f}
\label e.compression
$$ 
on $L^2(\T^d)$ viewed as a subspace of $L^2(\T^{d+1})$. For later use, let 
$$
H_f:=\{g \in L^2(\T^{d+1}) \st g = 0 \hbox{ a.e.\ on } (A_f)^c \}
$$
be the image of $T_{A_f}$.

We next remind the reader of the notion of stochastic domination between 
two probability measures
on $2^E$. First, if $\eta,\delta$ are elements of $2^E$, we write 
$\eta \preccurlyeq\delta$ if $\eta(e) \le \delta(e)$ for all $e\in E$. A 
subset $A$ of
$2^E$ is called {\bf increasing} if $\eta\in A$ and $\eta
\preccurlyeq\delta$ imply that $\delta\in A$.
If $\nu$ and $\mu$ are two probability measures on $2^E$, we write $\nu
\preccurlyeq \mu$ if $\nu(A) \le \mu(A)$ for all increasing sets $A$.
A theorem of \ref b.Strassen/ says that this is equivalent to the 
existence of a probability
measure $m$ on $2^E \times 2^E$ that has $\nu$ and $\mu$ as its first and
second marginals (i.e., $m$ is a {\bf coupling} of $\nu$ and $\mu$)
and such that $m$ is concentrated on the set
$\{(\eta,\delta):\eta \preccurlyeq\delta\}$ (i.e., $m$ is {\bf monotone}).

Throughout the paper, we shall use 
the following consequence of conditional negative association (\ref
t.negass/), sometimes called
{\bf joint negative regression dependence}
(see \ref b.Pemantle:negass/).

\procl p.JNRD
Assume that $\{X_i\}_{i\in I}$ has conditional negative association,
$I$ is the disjoint union of $I_1$ and $I_2$, and 
$a,b\in \{0,1\}^{I_1}$ with $a_i\le b_i$ for each $i\in I_1$.
Then 
$$
[\{X_i\}_{i\in I_2}\mid X_i=b_i, i\in I_1]
\preccurlyeq [\{X_i\}_{i\in I_2}\mid X_i=a_i, i\in I_1]\,,
$$
where $[ Y \mid A]$ stands for the law of $Y$ conditional on $A$.
\endprocl

\procl l.mon
Let $f_1,f_2:\Td\to [0,1]$ with $f_1\le f_2$ a.e. Then
$\P^{f_1} \preccurlyeq  \P^{f_2}$.
\endprocl 

\proof
This follows from a more general result (see \ref b.L:det/) that 
says that if $Q_1$ and $Q_2$ are two commuting positive contractions 
on $\ell^2(E)$ such that $Q_1\le Q_2$ in the sense that $Q_2-Q_1$ 
is positive, then $\P^{Q_1} \preccurlyeq  \P^{Q_2}$. 
However, here is a more concrete proof in our case. Since
$f_1\le f_2$ a.e., it follows that $H_{f_1}\subseteq H_{f_2}$, which implies
by Theorem 6.2 in \ref b.L:det/ that the projection measures 
$\P^{H_{f_1}}$ and $\P^{H_{f_2}}$ on $\Z^{d+1}$ satisfy 
$\P^{H_{f_1}}  \preccurlyeq  \P^{H_{f_2}}$ and hence their restrictions to $\Z^d$ satisfy
the same relationship; i.e., $\P^{f_1} \preccurlyeq\P^{f_2}$.
\Qed

We close this section with a key
representation of certain conditional probabilities and an application.
The minimum in \ref e.closest/ below is often referred to as a Szeg\H{o}
infimum.
For an infinite set $B \subseteq \Z^d \setminus \{ \bfz \}$, write
$$
\P^f[\eta(\bfz) = 1 \mid \eta \restrict B \equiv 1]
:=
\lim_{n \to\infty}
\P^f[\eta(\bfz) = 1 \mid \eta \restrict B_n \equiv 1]
\,,
$$
where $B_n$ is any increasing sequence of finite subsets of $B$ whose union is
$B$. This is a decreasing limit by virtue of \ref p.JNRD/.
Let $\ip{\cdot, \cdot}_f$ denote the usual inner product in the complex
Hilbert space $L^2(f)$.
For any set $B \subset \Z^d$, write $[B]_f$ for the closure in $L^2(f)$
of the linear span of the complex exponentials $ \{ \ee_k \st k
\in B \} $.

\procl t.closest
Let $f : \Td \to [0, 1]$ be measurable and $B\subset\Z^d$ with
$\bfz \notin B$.
Then 
$$
\P^f[\eta(\bfz) = 1 \mid \eta \restrict B \equiv 1]
=
\min \left\{ \int |\bfo - u|^2 f \,d\lambda_d \st u \in [B]_f \right\}
\,.
\label e.closest
$$
\endprocl

\proof
It suffices to prove the theorem for $B$ finite since the infinite case then
follows by a simple limiting argument.
So assume that $B$ is finite.
Note that $\widehat f(k-j) = \ip{\ee_j, \ee_k}_f$, so that 
$$
\P^f[\eta \restrict B \equiv 1] = 
\det [ \ip{\ee_j, \ee_k}_f]_{j, k \in B}
\,,
$$
and similarly for $\P^f[\eta \restrict (B \cup \{ \bfz \}) \equiv 1]$.
Thus, the left-hand side of \ref e.closest/ is a quotient of
determinants. The fact that such a quotient has the form of the right-hand
side is sometimes called Gram's formula. We include the proof
for the convenience of the reader.
Since $\ee_\bfz = \bfo$, it follows by row operations on the matrix
$[\ip{\ee_j, \ee_k}]_{j, k \in B \cup \{ \bfz \}}$ that
$$
\P^f[\eta(\bfz) = 1 \mid \eta \restrict B \equiv 1]
=
\|P_{[B]_f}^\perp \bfo\|_f^2
\,,
\label e.byproj
$$
where $P_{[B]_f}^\perp$ denotes orthogonal projection onto the orthogonal
complement of $[B]_f$ in $L^2(f)$.
Since this is the squared distance from $\bfo$ to $[B]_f$,
the equation \ref e.closest/ now follows.
\Qed

An extension of the above reasoning, given in \ref b.L:det/, provides
the entire conditional probability measure:

\procl t.condonB
Let $f : \Td \to [0, 1]$ be measurable and $B\subset\Z^d$.
Then the law of $\eta \restrict (\Zd \setminus B)$ conditioned on
$\eta \restrict B \equiv 1$ is the determinantal probability measure
corresponding to the positive contraction on $\ell^2(\Zd \setminus B)$
whose $(j, k)$-matrix entry is 
$$
\bigip{P_{[B]_f}^\perp \ee_j, P_{[B]_f}^\perp \ee_k}_f
\,.
$$
for $j, k \notin B$.
\endprocl

The Szeg\H{o} infimum that appears in \ref t.closest/ involves
trigonometric approximation, a classical area that
has strong connections to the topics of prediction and interpolation
for wide-sense stationary processes. We briefly discuss these topics
now. In later sections, we describe more explicit
connections to our results.
Recall that a mean-0 wide-sense stationary process is a
(not necessarily stationary) process
$\Seq{Y_n}_{n\in \Zd}$ for which all the variables have finite
variance and mean 0, and such that for each $k \in \Zd$, the covariance
$\Cov(Y_{n+k}, Y_n) = \E[Y_{n+k} \overline{Y_n}]$ does not depend on $n$.
There is then a positive measure $G$ (called the spectral measure) 
on $\T^d$ satisfying
$$
\widehat{G}(k)=\Cov(Y_{n+k}, Y_n)
$$
for $n,k\in \Zd$.
It turns out that for a one-dimensional wide-sense stationary process,
if $G$ is absolutely continuous with density $g$, then
$\GM(g)=0$ iff 
perfect linear prediction is possible, which means that $Y_0$ is in the closed
linear span of $\{Y_n\}_{n\le -1}$ in $L^2(\Omega)$,
where $\Omega$ is the underlying probability
space. This was proved in various versions by
Szeg\H{o}, Kolmogorov and Kre{\u\i}n.
Since the covariance with respect to $\P^f$ of 
$\eta(\bfz)$ and $\eta(k)$ for $k \in \Z^d$ is also
given by a Fourier coefficient (namely $-|\widehat f(k)|^2$ for $k \ne
\bfz$)
and since, as we shall see, the geometric mean of $f$ will play an 
important role in classifying the behavior of $\P^f$ as well,
one might wonder about the relationship between our
determinantal processes and wide-sense stationary processes.
It is not hard to show that $\P^f$ 
(viewed as a wide-sense stationary process)
has a spectral measure $G$ that is absolutely continuous with
density $g$ given by the formula
$$
g:= \widehat{f}(\bfz) \bfo -f*\widetilde{f}
\,,
$$
where $*$ denotes convolution and
$\widetilde{f}(t):=f(1-t)$.
Other than the trivial cases $f = \bfz$ and $f = \bfo$,
it is easy to check that $g$ is bounded away from 0 
and so, in particular, its geometric mean is always
strictly positive.
This suggests that our results are perhaps not so connected to
prediction and interpolation. However, it turns out that
some of the questions that
we deal with here (such as phase multiplicity and domination)
concerning $\P^f$ do have interpretations in terms of prediction and
interpolation for wide-sense stationary processes whose spectral
density is $f$ (which does not include $\P^f$).
More specifics will be given in the relevant sections.

Special attention is often devoted to stationary Gaussian processes, one
reason being that their distribution is determined entirely by their
spectral measure. It is known that
a stationary Gaussian process with no deterministic component
is a multistep Markov chain iff its spectral density is the reciprocal of a
trigonometric polynomial; see \ref b.Doob:elem/. 
The analogous property for determinantal processes is regeneration:

\procl p.regenerate
If $d=1$, then $f$ is the reciprocal of a trigonometric polynomial of
degree at most $n$ iff $\P^f$ is a regenerative process that regenerates
after $n$ successive $1$s appear.
\endprocl

This last property means that for any $k$, given that
$\eta \restrict [k+1, k+n] \equiv 1$, the future, $\eta \restrict
\CO{k+n+1,\infty}$, is conditionally independent of the past, $\eta
\restrict \OC{-\infty, k}$.

\proof
Note first that because of \ref t.condonB/, this
regenerative property holds for $\P^f$ iff for all $j \ge 0$ and all $C
\subset \OC{-\infty, -n-1}$, we have $P_{[B]_f}^\perp \ee_j = P_{[B \cup
C]_f}^\perp \ee_j$, where $B := [-n, -1]$.
(Here, we are relying on the fact that $\|P_{[B]_f}^\perp \ee_j\|_f > \|P_{[B
\cup C]_f}^\perp \ee_j\|_f$ if the vectors are not equal.)
This is the same as $P_{[B]_f}^\perp \ee_j \perp [C]_f$ for all $j \ge 0$, or,
in other words, there exists some $u_j \in [B]_f$ with $\ee_j - u_j \perp [B
\cup C]_f$. As this would have to hold for all $C$
(and $u_j$ is independent of $C$), 
it is the same as the existence of
some $u_j \in [B]_f$ such that $T_j := (\ee_j - u_j) f$ is analytic, i.e.,
$\widehat {T_j}(k) = 0$ for all $k < 0$. 
Now this implies that $\overline{T_0} = (\bfo - \overline{u_0}) f$, whence 
$T_0 (\bfo - \overline{u_0}) = \overline{T_0} (\bfo - u_0)$.
The left side of this last 
equation is an analytic function, while the right side
is the conjugate of an analytic function (i.e., its Fourier coefficients
vanish on $\Z^+$).
Therefore, both equal some constant, $c$.
Hence 
$$
f = {T_0 \over \bfo - u_0} = {c \over (\bfo - u_0) (\bfo - \overline{u_0})}
= {c \over |\bfo - u_0|^2}
\,,
$$
which is indeed the reciprocal of a trigonometric polynomial of degree at most
$n$.

Conversely, if $1/f$ is the reciprocal of a trigonometric polynomial of degree
at most $n$, then since $f \ge 0$, the theorem of Fej\'er and Riesz
(see \ref b.GrenanderSzego/, p.~20)
allows us to write $f = c/|\bfo - u_0|^2$ for some conjugate-analytic
polynomial $u_0 \in [B]_f$ such that the analytic extension of
$\bfo - \overline{u_0}$ to the unit disc has no zeroes.
We may rewrite this as $(\bfo - u_0) f = T_0$ for $T_0 := c/(\bfo -
\overline{u_0})$. 
Since $T_0$ has an extension to the unit disc as the reciprocal of an
analytic polynomial with no zeroes, it follows that $T_0$ is also analytic.
Multiplying both sides of this equation by $\ee_1$ and rewriting $\ee_1 u_0 f =
c_1 (\bfo - u'_1) f = c_1 (u_0 - u'_1) f + c_1 T_0$ for some constant $c_1$ and
some $u'_1 \in [B]_f$, we see that $(\ee_1 - u_1) f = T_1$ for $u_1 := c_1 (u_0
- u'_1) \in [B]_f$ and $T_1 := \ee_1 T_0 + c_1 T_0$, an analytic function.
Similarly, we may establish by induction that for each $j \ge 0$,
there is some $u_j \in [B]_f$ such that $T_j := (\ee_j - u_j) f$ is 
analytic.
This proves the equivalence desired.
\Qed

We now use this proof to establish the explicit probabilistic form of the
renewal process in \ref x.renewal/.
In this case, $B = \{-1\}$ and one easily verifies that $u_j = a^{j+1}
\ee_{-1}$, as is standard in the theory of linear prediction. 
Therefore, 
$$
\P^f[\eta(j) = 1 \mid \eta(-1) = 1]
=
\|P_{[B]_f}^\perp \ee_j\|_f^2
=
\|\ee_j - u_j\|_f^2
=
\|\ee_j\|_f^2 - \|u_j\|_f^2
=
{1-a \over 1+a} (1 - a^{2j+2})
\,.
$$
It now suffices to verify that this is also true for the explicit renewal
process described in \ref x.renewal/.
First, it is well known from basic renewal theory that for a renewal 
process, the probabilities
$\P^f[\eta(j) = 1 \mid \eta(-1) = 1]$ determine the distribution of the 
number of 0s between two 1s. Hence to verify the above statement,
one simply needs to check that these latter probabilities are related to
the interrenewal distribution via the appropriate convolution-type equation.
In this case, it comes down to verifying that for all $j\ge 0$,
$$
{1-a \over 1+a}(1-a^{2j+2})
=
(j+1)(1-a)^2a^j
+{1-a \over 1+a}\sum_{k=1}^j k (1-a)^2a^{k-1}(1-a^{2(j-k)+2})\,.
$$
This identity is easy to check.

\bsection {The Bernoulli Shift Property}{s.bern}

In this section,
we prove that the stationary determinantal processes studied here
are Bernoulli shifts. We 
assume the reader is familiar with the
basic notions of a Bernoulli shift (see, e.g., \ref b.Ornstein:book/).

\procl t.Bern
Let $f : \T^d \to [0, 1]$ be measurable.
Then $\P^f$ is a Bernoulli shift; i.e., it is isomorphic 
(in the sense of ergodic theory) to an i.i.d.\ process.
\endprocl

Before beginning the proof, we present a few preliminaries.
We first recall the definition of the $\dbar$-metric.

\procl d.overd
If $\mu$ and $\nu$ are 
$\Z^d$-invariant probability measures
on $2^{\Z^d}$, then
$$
\dbar(\mu,\nu) := \inf_m m \Big[ \big \{ (\eta,\delta) \in 2^{\Z^d} \times
2^{\Z^d} \st \eta(\bfz)\ne \delta(\bfz) \big \}\Big]\, ,
$$
where the infimum is taken over all couplings $m$ of $\mu$ and 
$\nu$ that are $\Z^d$-invariant.
\endprocl 

The following is a slight generalization of a well-known result
(see, for example, page 75 in \ref b.Liggett/). The 
proof is an immediate consequence of the existence of a monotone coupling
and is left to the reader.

\procl l.weak
Suppose that $\sigma_1$, $\sigma_2$, $\nu$ and $\mu$ are
$\Z^d$-invariant probability measures on $2^{\Z^d}$ such that $\nu
\preccurlyeq \sigma_i \preccurlyeq  \mu$ for both $i=1, 2$. 
Then 
$$
\dbar(\sigma_1, \sigma_2) \le \mu[\eta(\bfz) = 1] - \nu[\eta(\bfz)= 1]
\,.
$$
\endprocl 

\procl p.L1givesdbar
Let $f, g : \T^d \to [0, 1]$ be measurable.
Then 
$$
\dbar(\P^f, \P^g) \le \int_{\T^d} |f-g| \,d\lambda_d
\,.
$$
\endprocl 

\proof
Because of \ref l.mon/, we may apply \ref l.weak/ to $\sigma_1 := \P^f$,
$\sigma_2 := \P^g$, $\nu := \P^{f \wedge g}$, and $\mu := \P^{f \vee g}$.
We obtain that 
$$
\dbar(\P^f, \P^g)
\le
\int_{\T^d} f \vee g \,d\lambda_d - \int_{\T^d} f \wedge g \,d\lambda_d
=
\int_{\T^d} |f-g| \,d\lambda_d
\,.
\Qed
$$

\proofof t.Bern
The first step is to approximate $f$ by 
trigonometric polynomials. Let $K_r$ be the $r$th Fej\'er kernel for $\T$, 
$$
K_r := \sum_{|j| \le r} \Big(1-{|j|\over r+1}\Big) \ee_j 
\,,
$$
and define $K^d_r(x_1, \ldots, x_d) := \prod_{i=1}^d K_r(x_i)$.
It is well known 
that $K^d_r$ is a positive summability kernel, so that
if we define $g_r$ by
$$
g_r:= f * K^d_r\,,
$$
then $0\le g_r \le 1$ and $\lim_{r\to\infty} g_r =f$ a.e.\ and in $L^1(\Td)$. 

Next, since each $K^d_r$ is a trigonometric polynomial, so is each $g_r$.
{}From this and \ref l.fourier/, it is easy to see that there is a constant 
$C$ such that $\P^{g_r}(A\cap B) = \P^{g_r}(A)\P^{g_r}(B)$ 
if $A$ is of the form $\eta\equiv 1 \hbox{ on } S_1$ and 
$B$ is of the form $\eta\equiv 1 \hbox{ on } S_2$ 
with $S_1$ and $S_2$ being finite sets 
and having distance at least $C$ between them.
{}From this, a standard argument in probability (a $\pi$-$\lambda$
argument) shows that if $A$ is any event depending only on
locations $S_1$ and if $B$ is any event depending only on locations $S_2$
with $S_1$ and $S_2$ possibly infinite sets
having distance at least $C$ between them, then
$\P^{g_r}(A\cap B) = \P^{g_r}(A)\P^{g_r}(B)$. A process 
with this property is called a {\bf finitely dependent} process 
and this property implies it is a Bernoulli shift 
(e.g., the so-called very weak Bernoulli condition 
is immediately verified; see \ref b.Ornstein:book/ for this definition for
$d=1$ and, e.g., \ref b.Steif/ for general $d$).

Since the processes that are Bernoulli shifts 
are closed in the $\dbar$ metric (see
\ref b.Ornstein:book/) and \ref p.L1givesdbar/ tells us that
$\lim_{r\to\infty} \dbar(\P^{g_r},\P^f) = 0$, we 
conclude that $\P^f$ is a Bernoulli shift.
\Qed

\procl r.WB
An important property of 1-dimensional 
stationary processes is the weak Bernoulli
(WB) property. This is also referred to as ``$\beta$-mixing" and ``absolute
regularity" in the literature. Despite its name, it is known that WB is
strictly stronger than Bernoullicity. It is easy to check that ``aperiodic''
regenerative processes and finitely dependent processes are WB. Hence,
by our earlier results, if $f$ is a trigonometric polynomial or the inverse
of a trigonometric polynomial (or if $\bfo-f$ is the inverse
of a trigonometric polynomial), then $\P^f$ is WB. This is subsumed by the
independent
work of \ref b.ShiTak:II/, who showed that $\P^f$ is WB whenever
$\sum_{n \ge 1} n |\widehat f(n)|^2 <\infty$.
The precise class of $f$ for which $\P^f$ is WB is not known.
We also note that it follows from 
\ref b.Smor/ that if
$f$ is a trigonometric polynomial, then $\P^f$ is finitarily isomorphic
to an i.i.d.\ process.
\endprocl

\bsection {Support of the Measures}{s.fs}

In this section, we show that all of the probability measures $\P^f$ have
full support except in two degenerate cases. We begin with the following
lemma.

\procl l.lindep
Let $H$ be a closed subspace of $\hE$ and let $A$ and $B$ be finite disjoint
subsets of $E$. Then
$$
\P^H[\eta(e)=1 \hbox{ for } e\in A, \,\, \eta(e)=0 \hbox{ for } e\in B] >0
$$
if and only if 
$\{P_H(e)\}_{e\in A}\cup\{P_{H^\perp}(e)\}_{e\in B}$ is linearly independent.
\endprocl

\proof
In the special case that $A$ or $B$ is empty, the result follows from
\ref b.L:det/.
In general,
$\{P_H(e)\}_{e\in A}\cup\{P_{H^\perp}(e)\}_{e\in B}$ is linearly independent
if and only if
$\{P_H(e)\}_{e\in A}$ and $\{P_{H^\perp}(e)\}_{e\in B}$ are each linearly
independent sets. Also,
$\P^H[\eta \restrict A \equiv 1,\, \eta \restrict B \equiv 0] >0$
if and only if
$\P^H[\eta \restrict A \equiv 1] >0$ and
$\P^H[\eta \restrict B \equiv 0] >0$ since 
$$
\P^H[\eta \restrict A \equiv 1, \,\, \eta \restrict B \equiv 0] \ge  
\P^H[\eta \restrict A \equiv 1] \P^H[\eta \restrict B \equiv 0] 
$$ 
by the negative association property, \ref t.negass/.
Thus, the general case follows from the special case.
\Qed

\procl t.support
$\P^f$ has full support
for any function $f:\Td\to [0,1]$ other than $\bfz$ or $\bfo$.
\endprocl

\proof
Since marginals of probability measures with full support clearly
have full support, it suffices to prove the result for
$\P^H$ when $H$ is a closed $\Z^d$-invariant
subspace of $\hD$ other than
$\hD$ or 0. If we translate over to $L^2(\T^d)$, we see, by \ref l.lindep/,
that it suffices to show that
if $A\subseteq \Td$ with $0< \lambda_d(A) < 1$
and $n_1,\ldots,n_k, m_1,\ldots,m_\ell$ are
all distinct elements of $\Z^d$, then
$$
\{\ee_{n_j} \I A\}_{1\le j \le k} \cup
\{\ee_{m_r} \I {A^c}\}_{1\le r \le \ell} 
$$
is linearly independent.
Suppose that
$c_1,\ldots,c_k, d_1,\ldots,d_\ell$ are complex numbers such that
$$
\sum_{j=1}^k c_j \ee_{n_j} \I A +
\sum_{r=1}^\ell d_r \ee_{m_r}  \I {A^c} = \bfz 
\,.
$$
From this it follows that
$\sum_{j=1}^k c_j \ee_{n_j}=0$ a.e.\ on $A$. 
Since $\lambda_d(A) >0$,
we have that $\sum_{j=1}^k c_j \ee_{n_j}$ is 0 
on a set of positive measure.
It is well known that this implies that $c_1,\ldots,c_k$ vanish (the proof
uses induction on $d$ and Fubini's theorem). Similarly, $d_1,\ldots,d_\ell$
vanish.
\Qed

\bsection {Domination Properties}{s.dom}

In this section, we study the question of
which product measures are stochastically dominated by $\P^f$ and
which product measures stochastically dominate $\P^f$. 
For simplicity, we give our first results for $d=1$, and only afterwards
describe how these results extend to higher
dimensions. We also treat a different 
notion of ``full domination" at the end of the
section.
Recall that for $f : \Td \to [0, 1]$, we define 
$$
\GM(f) := \exp \int_\Td \log f \,d\lambda_d
\,.
$$

We introduce an auxiliary stronger
notion of domination than $\preccurlyeq$, but only in the case where
one of the measures is a product measure. Let $\mu_p$ denote product
measure with density $p$. 

\procl d.sdom
A stationary process $\{\eta_n\}_{n\in\Z}$ with distribution $\nu$
{\bf strongly dominates} $\mu_p$, written 
$\mu_p\strle\nu$,
if for any $n$ and any $a_1,\ldots,a_n \in \{0,1\}$,
$$
\P[\eta_0=1 \mid \eta_i=a_i, i=1,\ldots, n]\ge p
\,.
$$
Similarly, we define $\nu\strle\mu_p$ if the above inequality holds
when $\ge$ is replaced by $\le$.
\endprocl

The following lemma is easy and well known; it is sometimes referred to as 
Holley's lemma.

\procl l.holley
If $\mu_p\strle \nu$, then $\mu_p\preccurlyeq\nu$.
\endprocl 

The converse is not true, as we shall see later (\ref r.support/).
In light of \ref l.mon/, we have that
if $p \le f \le q$, then 
$$
\mu_p \preccurlyeq \P^f \preccurlyeq \mu_q
\,.
$$
The optimal improvement of these stochastic bounds is as follows.

\procl t.dom
For any $f : \T \to [0, 1]$, we have
$\mu_p\preccurlyeq\P^f$ iff $p\le \GM(f)$ iff $\mu_p\strle\P^f$. Similarly,
$\P^f\preccurlyeq\mu_q$ iff $q\ge 1-\GM(\bfo-f)$ iff $\P^f\strle\mu_q$.
In addition, for any stationary process $\mu$ that 
has conditional negative association, we have
$\mu_p\preccurlyeq\mu$ iff $\mu_p\strle\mu$.
\endprocl 

\proof
Let $d_n := D_{n-1}(f)$ be the probability of having $n$ 1s in a row.
According to Szeg\H{o}'s theorem (\ref b.GrenanderSzego/, pp.~44, 66),
$d_{n+1}/d_n$ is decreasing in $n$
and
$$
\lim_{n\to\infty}d_{n+1}/d_n= \GM(f)
= \lim_{n \to\infty} d_n^{1/n}
\,.
$$
In particular, $d_{n+1}/d_n\ge \GM(f)$ for all $n$.

\ref p.JNRD/ implies that
for any fixed $n$,
$$
\P^f[\eta_0=1 \mid \eta_i=a_i, i=1,\ldots, n]
$$
is minimized among all $a_1,\ldots,a_n \in \{0,1\}$
when $a_1 = a_2 = \cdots = a_n = 1$.
In this case, the value is $d_{n+1}/d_n$.
Since this is at least $\GM(f)$, 
we deduce that if $p\le \GM(f)$, then
$\mu_p\strle \P^f$ and hence that $\mu_p\preccurlyeq\P^f$. 

Conversely, if $\mu_p \preccurlyeq \P^f$, then certainly
$p^n\le d_n$ for all $n$. 
Hence $p \le \GM(f)$. The second to last statement can be proved in 
the same way or can be concluded by symmetry. 
Finally, the last statement can be proved in a similar fashion.
\Qed


\procl r.examples
The above theorem gives us two interesting examples. 
First, if we take $f:=\I A$ where $A$ has Lebesgue measure $1-\epsilon$, then
$\dbar(\P^f, \delta_{\bfo}) \le \epsilon$ by 
\ref p.L1givesdbar/, but nonetheless, $\P^f$ does not dominate
any nontrivial product measure since $\GM(f) = 0$.
Second, if we take
$f:={\bf 1}_{[0,1/2]}+ .4\cdot {\bf 1}_{[1/2,1]}$, 
then $f < 1/2$ on a set of positive
measure, but nonetheless $\P^f$ dominates $\P^{1/2}=\mu_{1/2}$
since $\GM(f) = (.4)^{1/2} > 1/2$.
\endprocl

\procl x.ustH_x-dom
Let $f$ be as in \ref x.ustH/, so that $\P^f$ is the law of the horizontal
edges of the uniform spanning tree in the plane.
In order to examine a 1-dimensional process, let us consider only
the edges lying on the $x$-axis.
If we let 
$$
g(x) := \int_\T f(x, y) \,d\lambda_1(y)
\,,
$$
then the edges lying on the $x$-axis have the law $\P^g$
since $\widehat g(k) = \widehat f(k, 0)$ for all $k \in \Z$.
Since an antiderivative of $1/(1 + a\sin^2 \pi y)$ is 
$$
{\arctan \left( \sqrt{1+a}\tan \pi y\right) \over \pi \sqrt{1+a}}
\,,
$$
we have 
$$
g(x) =
{\sin \pi x \over \sqrt{1+\sin^2 \pi x}}
\,,
$$
as given in \ref x.ustHx/.
We claim that $\P^g$ strongly dominates $\mu_p$ for $p := \sqrt 2 - 1$, and
this is optimal.
In order to show this, we calculate $\GM(g)$.
Write $g_1(x) := \sin^2 \pi x$. 
Then $(\GM(g))^2 = \GM(g_1)/\GM(\bfo + g_1)$.
Let $G_1$ be the harmonic extension of $\log g_1$ from the circle
to the unit disc.
Then $\int \log g_1 \,d\lambda_1 = G_1(0)$ by the mean value property of
harmonic functions.
Since $g_1 = |(1-\ee_1)/2|^2$, we see that $G_1$ is the
real part of the analytic function $z \mapsto 2 \log [(1-z)/2]$
in the unit disc, from which we conclude that $G_1(0) = \log (1/4)$.
Therefore, $\GM(g_1) = 1/4$. 
Similarly, $\bfo+g_1 = |[\sqrt 2 + 1 - (\sqrt 2 - 1)\ee_1]/2|^2$, whence
$\GM(\bfo+g_1) = [(\sqrt 2 + 1)/2]^2$.
Therefore $\GM(g) = \sqrt 2 - 1 = 0.4142^+$, as desired.
It turns out that $q := 1 - \GM(\bfo - g) = 1 - 2(\sqrt 2 -1)e^{-2\cat/\pi} =
0.5376^+$, where 
$$
\cat := \sum_{k=0}^\infty {(-1)^k \over (2k+1)^2} = 0.9160^-
\label e.defcat
$$
is {\bf Catalan's constant}.
Thus, $\P^g \strle \mu_q$.
It is interesting how close $p$ and $q$ are.
\endprocl

We now introduce some new mixing conditions.
Given a positive integer $r$, we may restrict $\P^f$ to $2^{r\Z}$.
If we identify $r\Z$ with $\Z$, then we obtain the process $\P^{f_r}$, where 
$$
f_r(t) := {1 \over r} \sum_{x \in r^{-1} t} f(x)
\,,
$$
where $r^{-1} t := \{ x \in \T \st r x = t \}$.
The reason that this restriction is equal to $\P^{f_r}$ is that
for all $k \in \Z$, we have $\widehat{f_r}(k) = \widehat f(rk)$, as
is easy to check.
Because of this relation, we have that $\int |f_r(t) - \widehat
f(0)|^2\,d\lambda_1(t) \to 0$ as $r \to\infty$, so that by \ref l.weak/, it
follows that $\dbar(\P^{f_r}, \P^{\widehat f(0)}) \to 0$ as $r \to\infty$.
(It is not hard to show that a similar property holds for any Kolmogorov
automorphism (see \ref d.K/), while the example in \ref r.support/ shows
that this property can occur in other cases as well.)
In fact, we often have a stronger convergence for our determinantal processes,
as we show next.

\procl t.dom-mix
Let $f : \T \to [0, 1]$ be measurable.
If $\GM(f) > 0$ or $f$ is bounded away from 0 on an interval of positive length,
then there exist constants $p_r \to \widehat f(0)$ such that $\P^{f_r}
\succcurlyeq \P^{p_r}$ for all $r$.
Therefore, if $\GM\big(f(\bfo - f)\big) > 0$ or if $f$ is continuous and not
equal to $\bfz$ nor $\bfo$, then there exist constants $p_r, q_r \to \widehat
f(0)$ such that $\P^{p_r} \preccurlyeq \P^{f_r} \preccurlyeq \P^{q_r}$ for all
$r$.
\endprocl

Of course, in view of \ref t.dom/, we use $p_r := \GM(f_r)$ and $q_r := 1 -
\GM(\bfo - f_r)$.
Thus, the theorem is an immediate consequence of the following lemma.

\procl l.dom-mix
Let $f : \T \to [0, 1]$ be measurable.
If $\GM(f) > 0$ or $f$ is bounded away from 0 on an interval of positive length,
then $\GM(f_r) \to \widehat f(0)$ as $r \to\infty$.
\endprocl

\proof
We have seen that $f_r \to \fav \cdot \bfo$ in $L^2$, whence also in measure.
Therefore $\log f_r \to (\log \fav) \cdot \bfo$ in measure.
Thus, it remains to show that $\{ \log f_r \}$ is uniformly integrable (at
least, for all large $r$).
Suppose first that $\GM(f) > 0$.

Given $h \in L^1(\T)$, write $h^{(r)}(t) := h(r t)$, so that
$\widehat{h^{(r)}}(k)$ is 0 when $k$ is not a multiple of $r$ and is
$\widehat h(k/r)$ when $k$ is a multiple of $r$.
For any $g, h \in L^2(\T)$, we have
\begineqalno
\int_\T g_r(t) \overline{h(t)} \,d\lambda_1(t)
&=
\sum_{k \in \Z} \widehat{g_r}(k) \overline{\widehat h(k)}
=
\sum_{k \in \Z} \widehat{g}(r k) \overline{\widehat h(k)}
=
\sum_{k \in \Z} \widehat{g}(k) \overline{\widehat{h^{(r)}}(k)}
\cr&=
\int_\T g(t) \overline{h^{(r)}(t)} \,d\lambda_1(t)
=
\int_\T g(t) \overline{h(r t)} \,d\lambda_1(t)
\,.
\label e.turn
\cr
\endeqalno
Therefore, for any set $A \subseteq \T$, we have 
$$
\int_{A} \log (f_r) \,d\lambda_1
\ge
\int_{A} (\log f)_r \,d\lambda_1
=
\int_{r^{-1} A} \log f\,d\lambda_1
\,;
$$
the inequality follows by concavity of log, while the equality follows from
\ref e.turn/ applied to $g := \log f$ and $h := \I A$.
Given any $\epsilon \in \big(0, \fav\big)$, let $A^r_\epsilon := \{ t \st
f_r(t) < \epsilon \}$.
Note that for any measurable $A \subset \T$, we have
$\lambda_1(r^{-1} A) = \lambda_1(A)$.
Since $f_r \to \fav \cdot \bfo$ in measure, it follows that
$$
\lim_{r \to\infty} \lambda_1(r^{-1} A^r_\epsilon) =
\lim_{r \to\infty} \lambda_1(A^r_\epsilon)
= 0
\,.
$$
Since $\log f$ is integrable, it follows that 
$$
\lim_{r \to\infty} \int_{r^{-1} A^r_\epsilon} \log f \,d\lambda_1 =
0
\,.
$$
This establishes the uniform integrability in the first case.

In the second case where $f \ge c > 0$ on an interval of length
$\epsilon > 0$, we have that $f_r \ge c \epsilon$ on all of $\T$ for
all $r > 1/\epsilon$.
It is then obvious that $\{\log f_r\}$ is uniformly integrable for all 
$r > 1/\epsilon$.
\Qed

\procl r.nasty
There are functions $f$ with $f > 0$ a.e., yet $\GM(f_r) = 0$ for all $r$.
For example, enumerate the rationals $\{ x_j \st j \ge 1\}$ in $(0, 1)$ and
choose $\epsilon_j > 0$ such that $x_j + \epsilon_j < 1$ and $\sum_j
\epsilon_j < 1$.
Define $A_n := [0, 1] \setminus \bigcup_{j > n} [x_j, x_j + \epsilon_j]$ and
$f_0(x) := e^{-1/x}$.
Then one can show that 
$$
f(x) := f_0(x) \I{A_0}(x) +
\sum_{n=1}^\infty f_0(x-x_n) \I{[x_n, x_n+\epsilon_n] \cap A_n}(x) 
$$
is such an example.
\endprocl

We turn next to the extension of the prior results to higher dimensions.
This is basically straightforward, but has interesting applications, as we
shall see.
First, we recall the usual notion of the {\bf past $\sigma$-field} $\past(k)$
at a point $k \in \Z^d$.
We use lexicographic ordering on $\Z^d$, i.e., write $(k_1, k_2,
\ldots, k_d) \prec (l_1, l_2, \ldots, l_d)$ if $k_i < l_i$ when $i$ is the
smallest index such that $k_i \ne l_i$.
Then define $\past(k)$ to be the $\sigma$-field generated by $\eta(j)$ for all
$j \prec k$.
More generally, the past $\sigma$-field could be defined with respect to 
any {\bf ordering} of $\Z^d$, which means the selection of a set $\ord \subset
\Z^d$ that has the properties $\ord \cup (-\ord) = \Zd \setminus \{ \bfz \}$, 
$\ord \cap (-\ord) = \emptyset$, and $\ord+\ord \subset \ord$. The associated ordering 
is that where $k \prec l$ iff $ l-k \in \ord$.
$\past_\ord(u)$ will denote the past of $u$ with respect to the ordering
$\ord$,
i.e., $\{v \st v \prec u\}$.
Thus, $\past_\ord(\bfz)$ is just $-\ord$. 
(For a characterization of
all orders, see \ref b.Teh/, \ref b.Zaiceva/, or \ref b.Trevisan/.)
As before, let $\mu_p$ denote product measure with density $p$.

\procl d.sdom-d
Given an ordering $\ord$, 
a stationary process $\{\eta_n\}_{n\in\Z^d}$ with distribution $\nu$
{\bf strongly dominates} $\mu_p$, written 
$\mu_p\strle\nu$,
if 
$$
\P[\eta_\bfz=1 \mid \past_\ord(\bfz)]\ge p \quad\nu\hbox{-a.s.}
$$
Similarly, we define $\nu\strle\mu_p$ if the above inequality holds
when $\ge$ is replaced by $\le$. (Note that the ordering $\ord$ here is 
suppressed in the notation.)
\endprocl

Although we have phrased it differently, in one dimension, this is equivalent
to \ref d.sdom/ when $\ord=\{-1,-2,\ldots\}$.
Again, we have a version of Holley's lemma:

\procl l.holley-d
Given any ordering $\ord$, 
if $\mu_p\strle \nu$, then $\mu_p\preccurlyeq\nu$.
\endprocl 

(Note that to produce a monotone coupling for a general ordering with
respect to a set $\ord$, it is enough to do so for the measures
restricted to any finite $B \subset \Zd$. Given such a $B$, order $B$ by the
restriction of $\prec$ to $B$ and couple the measures by adding sites from
$B$ in this order.)

We now prove

\procl t.dom-d
Fix any ordering $\ord$.
For any measurable $f : \Td \to [0, 1]$, we have
$\mu_p\preccurlyeq\P^f$ iff $p\le \GM(f)$ iff $\mu_p\strle\P^f$. Similarly,
$\P^f\preccurlyeq\mu_q$ iff $q\ge 1-\GM(\bfo-f)$ iff $\P^f\strle\mu_q$.
In addition, for any stationary process $\mu$ that
has conditional negative association, we have
$\mu_p\preccurlyeq\mu$ iff $\mu_p\strle\mu$.
\endprocl 

\proof
\ref p.JNRD/ and \ref t.closest/ imply that 
\begineqalno
\essinf \P^f\Big[\eta(\bfz) = 1 \mid \past_\ord(\bfz)\Big] 
&=
\inf \Big\{ \P^f[\eta(\bfz) = 1 \mid \eta \restrict A \equiv 1] \st
A \subset -\ord \hbox{ is finite} \Big\}
\cr&=
\P^f[\eta(\bfz) = 1 \mid \eta \restrict (-\ord) \equiv 1]
\cr&=
\min \left\{ \int |\bfo - T|^2 f \,d\lambda_d \st T \in [-\ord]_f \right\}
\,.
\endeqalno
By \ref b.HelLow:I/, the latter quantity equals $\GM(f)$.
Therefore, $p\le \GM(f)$ iff $\mu_p\strle\P^f$. 
On the other hand, if $\mu_p\preccurlyeq\P^f$, then let $A(r)$ be the
finite subset of $-\ord$ consisting of all points within some large radius $r$
about $\bfz$.
Let $a_1 \prec a_2 \prec \cdots \prec a_n$ be the elements of $A(r)$ in
order, where $n := |A(r)|$, and write $A_j := A(r) \cap (-\ord+a_j)$.
Then 
$$
p^n \le \P^f[\eta \restrict A(r) \equiv 1]
=
\prod_{j=1}^n \P^f[\eta(a_j) = 1 \mid \eta \restrict A_j \equiv 1]
\,.
$$
Most terms in this product are quite close to $\GM(f)$, while all lie in
$[\GM(f), 1]$.
It follows that
$$
\lim_{r \to\infty} \P^f[\eta \restrict A(r) \equiv 1]^{1/n}
=
\GM(f)
\,,
$$
so that $p \le \GM(f)$.
Finally, the remaining statements can be proved in a similar fashion.
\Qed

\procl r.past
It follows from \ref t.dom-d/ that the choice of ordering $\ord$ does not
determine whether $\mu_p \strle \P^f$.
This is not true for general stationary processes, even in one dimension. 
We are grateful to Olle H\"aggstr\"om for the following simple example.
Let $\mu$ be the distribution on $\{0, 1 \}^2$ given by 
$$
\mu = (1/7) (\delta_{(0,0)} +\delta_{(0,1)}) + (2/7)\delta_{(1,0)} + (3/7)
\delta_{(1,1)}
\,.
$$
Let $\eta \in \{ 0, 1 \}^\Z$ be such that with probability 1/2, $(\eta_{2n},
\eta_{2n+1})$ are chosen independently each with distribution $\mu$ and
otherwise $(\eta_{2n-1}, \eta_{2n})$ are chosen independently each with
distribution $\mu$.
Let $\nu$ be the law of $\eta$.
If $\ord := \Z^-$, then $\mu_p \strle \nu$ iff $p \le 1/2$, while
if $\ord := \Z^+$, then $\mu_p \strle \nu$ iff $p \le 4/7$.
\endprocl

\procl x.ustH-dom
Let $f$ be as in \ref x.ustH/, so that $\P^f$ is the law of the horizontal
edges of a uniform spanning tree in the square lattice $\Z^2$.
By \ref b.Kas:dim/ or \ref b.Montroll/, we have
$$
\int_{\T^2} \log 4(\sin^2 \pi x + \sin^2 \pi y) \,d\lambda_2(x, y)
=
{4\cat \over \pi}
=
1.1662^+
\,,
\label e.catalan
$$
where $\cat$ is again Catalan's constant \ref e.defcat/.
As we have shown in \ref x.ustH_x-dom/, $\GM(\sin^2 \pi x)=1/4$,
whence by \ref e.catalan/,
$\GM(f) = e^{-4\cat/\pi} =\GM(\bfo - f)$, 
where we are using the observation that $1-f(x, y) = f(y,x)$. 
(This last identity has a combinatorial reason arising from planar dual
trees.) 
Therefore $\mu_p \preccurlyeq \P^f\preccurlyeq \mu_{1-p}$ for $p :=
e^{-4\cat/\pi} = 0.3115^+$ and this $p$ is optimal (on each side).
This result in itself is rather surprising and it would be fascinating to see
an explicit monotone coupling.
\endprocl

For $f:\Td\to [0,1]$ measurable and
$r = (r_1, r_2, \ldots, r_d) \in \Z^d$ with all $r_j > 0$, define
$$
f_r(t) := {1 \over \prod_{j=1}^d r_j} \sum_{x \in r^{-1} t} f(x)
\,,
$$
where $r^{-1} (t_1, t_2, \ldots, t_d) := \{ (x_1, x_2, \ldots, x_d) \in \Td
\st \all j\ r_j x_j = t_j \}$.
Write $r \to\infty$ to mean that $\min r_j \to\infty$.
The following is a straightforward extension of \ref t.dom-mix/.
We leave its proof to the reader.

\procl t.dom-mix-d
Let $f : \Td \to [0, 1]$ be measurable.
If $\GM(f) > 0$ or $f$ is positive on a non-empty open set,
then there exist constants $p_r \to \widehat f(\bfz)$ as $r \to\infty$
such that $\P^{f_r}
\succcurlyeq \P^{p_r}$ for all $r$.
Therefore, if $\GM\big(f(\bfo - f)\big) > 0$ or if $f$ is continuous and not
equal to $\bfz$ nor $\bfo$, then there exist constants $p_r, q_r \to \widehat
f(\bfz)$ such that $\P^{p_r} \preccurlyeq \P^{f_r} \preccurlyeq \P^{q_r}$ for
all $r$.
\endprocl

Consider now a domination property even stronger than our previously defined
strong domination. 
Let $\others$ be the $\sigma$-field generated by $\eta(k)$ for $k \ne \bfz$.

\procl d.fulldom
A stationary process $\{\eta_n\}_{n\in\Z^d}$ with distribution $\nu$
{\bf fully dominates} $\mu_p$, written 
$\mu_p\fullle\nu$,
if 
$$
\P[\eta_\bfz=1 \mid \others]\ge p \quad\nu\hbox{-a.s.}
$$
In this situation, we also say that $\nu$ is {\bf uniformly insertion tolerant
at level} $p$.
Similarly, we define $\nu\fullle\mu_p$ if the above inequality holds
when $\ge$ is replaced by $\le$, and
say that $\nu$ is {\bf uniformly deletion tolerant at level} $1-p$.
We say that $\nu$ is {\bf uniformly insertion tolerant} if
$\mu_p\fullle\nu$ for some $p > 0$ and
that $\nu$ is {\bf uniformly deletion tolerant}
if $\nu\fullle\mu_p$ for some $p < 1$.
\endprocl

We show that the optimal level of uniform insertion tolerance of a
determinantal process $\P^f$
is the {\bf harmonic mean} of $f$, defined as 
$$
\HM(f) := \left(\int_\Td {d\lambda_d  \over f} \right)^{-1}
\,.
$$
Note that $\HM(f) = 0$ iff $1/f$ is not integrable.

\procl t.fulldom
For any measurable $f : \Td \to [0, 1]$, we have
$\mu_p\fullle\P^f$ iff $p\le \HM(f)$. Similarly,
$\P^f\fullle\mu_q$ iff $q\ge 1-\HM(\bfo-f)$.
\endprocl

\proof
By \ref p.JNRD/ and \ref t.closest/, we have
\begineqalno
\essinf \P^f\Big[\eta(\bfz) = 1 \mid \others\Big] 
&=
\inf \Big\{ \P^f[\eta(\bfz) = 1 \mid \eta \restrict A \equiv 1] \st
A \subset \Zd \setminus \{ \bfz \} \hbox{ is finite } \Big\}
\cr&=
\P^f[\eta(\bfz) = 1 \mid \eta \restrict B \equiv 1]
=
\min \left\{ \int |\bfo - T|^2 f \,d\lambda_d \st T \in [B]_f \right\}
\, ,
\endeqalno
where $B := \Zd \setminus \{ \bfz \}$.
As shown by \refbmulti{Kolmog:sshs,Kolmog:iesrs} for $d=1$, the latter equals
$\HM(f)$.
The proof extends immediately to general $d$.
This proves the first assertion. The second follows by symmetry.
\Qed

\procl r.HM
Since a proof of Kolmogorov's theorem that 
$$
\min \left\{ \int |\bfo - T|^2 f \,d\lambda_d \st T \in [B]_f \right\}
=
\HM(f)
\,,
$$
where $B := \Zd \setminus \{ \bfz \}$,
is difficult to find in readily accessible sources, we
provide one here.
We have
$$
\min \left\{ \int |\bfo - T|^2 f \,d\lambda_d \st T \in [B]_f \right\}
=
\|u\|^2_f\,,
$$
where
$$
u :=
\|P^\perp_{[B]_f} \bfo\|^2_f
\,.
$$
Now $g \perp [B]_f$ iff $g \in L^2(f)$ and $\widehat{gf}(k) = 0$ for all $k
\in B$.
The latter condition holds iff $gf$ is a constant.
If $g \ne \bfz$, then we deduce that $1/f \in L^2(f)$, i.e., $\HM(f) > 0$.
Therefore,
$$
[B]^\perp_f = \cases{0 &if $\HM(f) = 0$,\cr
               \C /f &if $\HM(f) > 0$.\cr
               }
$$
Hence $u = \bfz$ if $\HM(f) = 0$, while otherwise, 
$$
\|u\|^2_f = |\ip{\bfo, \sqrt{\HM(f)}/f}_f|^2
=
\HM(f)
\,,
$$
as desired.
\endprocl

\procl r.strongnotfull
It is easy to find $f$ such that $\GM(f)> 0$ and $\HM(f)= 0$.
Indeed, a natural such example is the function $g$ of \ref x.ustHx/.
For any such $f$, the corresponding process
$\P^f$ strongly dominates a
nontrivial product measure, but does not
fully dominate any nontrivial product measure. In addition, for any function
$f$ such that $\HM(f) > 0$
and $f$ is not constant a.e., $\GM(f)> \HM(f)$ (as a consequence of
Jensen's inequality), so that $\P^f$ will strongly dominate strictly more
product measures than it will fully dominate.
\endprocl

\procl x.ustH-delete
Let $f$ be as in \ref e.ustH-d/, so that $\P^f$ is the law of the edges of the
uniform spanning forest in $\Z^d$ that lie parallel to the $x_1$-axis.
Then $\P^f$ is uniformly deletion tolerant 
iff $d \ge 4$.
This is because $1/(\bfo - f)$ is integrable iff $d \ge 4$.
For example, when $d=4$, we obtain full domination by $\mu_p$ with $p :=
0.66425^-$, where we have calculated $\HM(\bfo-f)$ as follows. First, we have
$$
1/(\bfo - f) = 
1 + {\sin^2 \pi x_1 \over \sum_{j=2}^4 \sin^2 \pi x_j}
\,,
$$
whence 
\begineqalno
\int_{\T^4} 1/(\bfo - f) \,d\lambda_4
&=
1 + {1 \over 2} \int_{\T^3} {1 \over \sum_{j=2}^4 \sin^2 \pi x_j}
\,d\lambda_3(x_2, x_3, x_4)
\cr&=
1 + {1 \over 2} \int_{\T^2} {1 \over \sqrt{\left(\sum_{j=2}^3 \sin^2 \pi x_j
\right) \left( 1 + \sum_{j=2}^3 \sin^2 \pi x_j\right)}} \,d\lambda_2(x_2, x_3)
\,.
\cr
\endeqalno
This last integral has no simple form and is calculated numerically.
This gives the value reported for $1 - \HM(\bfo - f)$.
For large $d$, we have $1 - \HM(\bfo - f) \sim 1/d$. Indeed, write
$\lambda_{d-1}^*$ for Lebesgue measure on $\R^{d-1}$. Letting
$$
A_d:=\int_{\T^{d-1}} {d - 1 \over \sum_{j=2}^{d} 2 \sin^2 \pi x_j}
\,d\lambda_{d-1}(x_2, \ldots, x_d)\,,
$$
we have that 
$$
(d-1) (1 - \HM(\bfo - f))= {A_d\over 1+A_d/(d-1)}
\,.
$$
Hence it suffices to show that $\lim_{d\to\infty}A_d=1$.
To show this, note that
\begineqalno
A_d
&=
\int_0^\infty \lambda_{d-1}\Big[ \sum_{j=2}^{d} 2 \sin^2 \pi x_j <
(d-1)/t\Big]
\,dt
\cr&\le
\int_0^2 \lambda_{d-1}\Big[ \sum_{j=2}^{d} 2 \sin^2 \pi x_j < (d-1)/t\Big] \,dt
\cr&\quad+ \int_2^{10(d-1)} 
\lambda_{d-1}\Big[ \sum_{j=2}^{d} 2 \sin^2 \pi x_j < (d-1)/2\Big] \,dt
\cr&\quad+ \int_{10(d-1)}^\infty
\lambda_{d-1}^*\Big[ \all j |x_j| < 1/2 \hbox{ and }
\sum_{j=2}^{d} 8 x_j^2 < (d-1)/t\Big] \,dt
\cr&\to
\int_0^2 \I{[0, 1]}(t) \,dt
=
1
\endeqalno
as $d \to \infty$, where we have used the weak law of large
numbers and the bounded convergence theorem for the first piece,
a standard large-deviation result for the second piece, and an easy 
estimate on the third piece. The reverse inequality obtains by using only
the first piece. Thus, $A_d \to 1$, as desired.
This value of $1 - \HM(\bfo - f)$,
which gives a full domination upper bound on $\P^f$, should be
compared to $\widehat f(\bfz) = 1/d$ (by symmetry), which is the
probability of a 1 at a site.
On the other hand, for no $d$ is the process uniformly insertion tolerant
since $1/f$ is never integrable.
We remark that for the full uniform spanning forest measure on $\Zd$
(considering edges in all directions), we have change intolerance, meaning
that $\P[\eta(\bfz) = 1 \mid \others] \in \{ 0, 1 \}$ a.s. This follows from a
result of \BLPSusf, as explained by \ref b.HeicklenLyons/.
\endprocl

\procl x.ust-delete
Let $g(x)$ be as in \ref x.ustHx/, so that
the edges of the uniform spanning tree in the plane that
lie on the $x$-axis have the law $\P^g$.
We have $\P^g \fullle \mu_p$ for $p := (1+\pi)/(1+2\pi) = 0.56865^+$, and this
is optimal.
This is because 
\begineqalno
\int_\T 1/(1-g) \,d\lambda 
&=
\int_\T \left(1 + \sin^2 \pi x  + \sin \pi x \sqrt{1 + \sin^2 \pi x }\right)
\,d\lambda(x)
\cr&=
3/2 + \int_\T \sin \pi x \sqrt{1 + \sin^2 \pi x } \,d\lambda(x)
\,.
\cr
\endeqalno
An antiderivative of this integrand is 
$$
-{1 \over \pi} \arctan\left({\cos\pi x \over \sqrt{1 + \sin^2\pi x}}\right)
-{1 \over 2\pi} \cos\pi x\sqrt{1 + \sin^2\pi x}
\,,
$$
whence the remaining integral is $1/2+1/\pi$. Therefore $\int_\T 1/(1-g)
\,d\lambda = 2+1/\pi$ and $1 - \HM(\bfo-g) = (1+\pi)/(1+2\pi)$, as desired.
Observe that $\int 1/g =\infty$ and so the process is not uniformly
insertion tolerant.
Similarly, for the edges lying on the $x$-axis of the uniform spanning tree in
3 dimensions, we have full domination by $\mu_p$ with $p := 0.37732^+$.
Here, the law of these edges is $\P^h$, where 
$$
h(x) :=
\int_{\T^2} f(x, y, z) \,d\lambda_2(y, z)
=
\int_\T {\sin^2 \pi x \over \sqrt{\left(\sin^2 \pi x + \sin^2 \pi y\right)
\left(1 + \sin^2 \pi x + \sin^2 \pi y\right)}} \,d\lambda_1(y)
\,.
$$
This integral has no simpler form, so to compute $p := 1 - \HM(\bfo - h)$, we
calculated $h$ numerically and used the result to calculate the harmonic mean
numerically. Since $\int 1/h =\infty$, as is easily checked,
this process is not uniformly insertion tolerant.
Similarly, one can check that the process of edges lying on the $x$-axis of the 
uniform spanning tree in $d\ge 4$ dimensions
is not uniformly insertion tolerant.
\endprocl

\procl x.USTzz-dom
Let $f$ be as in \ref x.USTzz/.
It turns out that $\GM(f) = e^{-2\cat/\pi}/\sqrt2 = 1 - \GM(\bfo-f) =   
0.39467^+$. This gives strong domination inequalities.
It is easy to see that $\HM(f) = \HM(\bfo - f) = 0$, so there are no
nontrivial full domination inequalities.
The interest of this function is that it describes the process of edges
along a zig-zag path in the plane for the uniform spanning tree, as claimed
in \ref x.USTzz/. 
We now sketch how to prove this.
The Transfer Current Theorem of \ref b.BurPem/ allows one to calculate, via
determinants, the law for any set of possible edges of the uniform spanning
tree.
Thus, it is enough to verify that for the edges belonging to the zig-zag
path, the matrix entries are those of the Toeplitz matrix associated to
the Fourier coefficients given in \ref x.USTzz/.
For even $k$, the values of $\widehat f(k)$ are given in
\ref b.Lyons:book/.
(These values imply the astonishing
fact that the edges in the plane that lie along a diagonal, e.g., the horizontal
edges with left endpoints $(n, n)$ ($n \in \Z$), are independent, i.e., have
law $\mu_{1/2}$.)
Thus, it remains to treat the case of odd $k$.
A straightforward application of the Transfer Current Theorem gives that
the matrix entry corresponding to the edges $e_0$ and $e_{2k-1}$ is 
\begineqalno
&\int {e(s) + e(t) - e(s+t) - 1 \over 4 - e(s) - e(-s) - e(t) - e(-t)}
e(k s + k t) \,d\lambda_2(s, t)
\cr&\qquad\qquad=
\int {e(s) + e(x-s) - e(x) - 1 \over 4 - e(s) - e(-s) - e(x-s) - e(s-x)}
e(k x) \,d\lambda_2(s, x)
\,,
\endeqalno
where $e(x) := e^{2 \pi i x}$.
Evaluate the integral in $s$ for fixed $x$ by a contour integral 
$$
{1 \over 2 \pi i} \oint {z + e(x) z^{-1} - e(x) - 1 \over 4 - z - z^{-1} -
e(x) z^{-1} - e(-x) z} {d z \over z}
$$
over the contour $|z| = 1$.
The integrand has poles inside the unit disc at $z = 0$ and 
$$
z = {2-\sqrt{4-|1+e(x)|^2} \over 1+e(-x)}
\,.
$$
After use of the residue theorem and integrating in $x$, one obtains
$\widehat f(2k-1)$, as desired.
\endprocl

We close this section by
describing how our domination results can be interpreted in terms of
prediction and interpolation questions for wide-sense stationary 
processes. In view of \ref t.closest/ and the well-known correspondence
between prediction and Szeg\H{o} infima, for $d=1$,
$\P^f[\eta(\bfz) = 1 \mid \eta \restrict \{-1,-2,\dots\} \equiv 1]$
is exactly the 
mean squared error for the 
best linear predictor of $Y_0$ given $\Seq{Y_n}_{n\le -1}$,
where 
$\Seq{Y_n}$ is a wide-sense stationary process with spectral density $f$.
Similarly,
$\P^f[\eta(\bfz) = 0 \mid \eta \restrict \{-1,-2,\dots\} \equiv 0]$
is exactly the mean squared error for the best linear predictor
of $Y_0$ given $\Seq{Y_n}_{n\le -1}$, where 
$\Seq{Y_n}$ is now a wide-sense stationary process with spectral density 
$\bfo-f$. This gives us a correspondence between strong domination and 
prediction. An analogous correspondence holds between
full domination and interpolation, where one instead looks at the 
mean squared error for the 
best linear predictor of $Y_0$ given $\Seq{Y_n}_{n\neq 0}$.  

\bsection{Entropy}{s.ent}

We assume the reader is familiar with the definition of the entropy
$H(\mu)$ of a process $\mu$, as well as basic results concerning entropy (see
\ref b.Walters/ and \ref b.KatzWeiss/).  
Because of Ornstein's theorem (and its generalizations, see \ref
b.KatzWeiss/, \ref b.Conze/, \ref b.Thouvenot/, and \ref b.OrnWeiss/)
that entropy characterizes
Bernoulli shifts up to isomorphism, the following question is particularly
interesting:

\procl q.entropy
What is $H(\P^f)$?
\endprocl

We know the answer only in the trivial case where $f$ is a constant and in
case $f$ or $\bfo - f$ is the reciprocal of a trigonometric polynomial of
degree 1. In principle, as we shall see, one can also determine the entropy
when $f$ or $\bfo - f$ is the reciprocal of any trigonometric polynomial,
but the formula would be rather unwieldy.

In general, then, we shall discuss how to estimate the entropy of $\P^f$.
The definition of entropy always provides upper bounds, due to
subadditivity, so the harder bound is the lower bound.
As we shall see, reciprocals of trigonometric polynomials can be used to
get arbitrarily close lower (and upper) bounds.
Unfortunately, that method is not practical for precise computation.
Nevertheless, in many cases, we have
another method that appears to work quite well.
Indeed, our method seems to provide
upper bounds that converge more quickly than does use of the definition.

\ref b.ShiTak:II/ proved that $H(\P^f) > 0$ for all $f \ne \bfz, \bfo$.
Of course, this also follows immediately from our \ref t.Bern/.

Let
$$
H[p]:= H(\mu_p)= -p\log p-(1-p)\log(1-p)
\,.
$$
\ref t.dom-d/ yields easy  
lower bounds on the entropy of those processes $\P^f$ such that $f(\bfo-f)$
has a strictly positive geometric mean. We shall obtain more refined lower
bounds later in the one-dimensional case. It is easy to see that
$$
\mu_p \strle \mu \strle \mu_{1-p}
\implies
H(\mu)\ge H[p]
\,.
\label e.doment
$$
By \ref t.dom-d/ and \ref e.doment/, we deduce the
following lower bound on entropy:

\procl p.entropy-d 
For any measurable $f : \Td \to [0, 1]$, we have
$$
H(\P^f)\ge \min\Big\{H\big[\GM(f)\big],H\big[\GM(\bfo-f)\big]\Big\}
\,.
$$
\endprocl 

\procl r.AMGM
This can be compared to the trivial upper bound $H(\P^f) \le H\big[\widehat
f(0)\big]$. Note that $\widehat f(0)$ is the arithmetic mean of $f$.
\ref b.ShiTak:II/ also note this upper bound and provide two lower bounds:
$H(\P^f) \ge H\big[\widehat f(0)\big]/2$ and
$H(\P^f) \ge \int_{\Td} H[f(x)] \,d\lambda_d(x)$.
\endprocl

\procl r.support
It is interesting that \ref e.doment/
is not true if $\strle$ is replaced by $\preccurlyeq$.
For example, let $(X_{2i},X_{2i+1})$ be, independently for 
different $i$, $(1,1)$ or $(0,0)$, each with probability $1/2$.
If $\nu$ is the distribution of this process, then $\nu$ is not
stationary, but $\mu:=\big(\nu +T(\nu)\big)/2$ is stationary, where $T$ is the
shift.
Next, $H(\mu)=\log 2/2 < H[1/\sqrt{2}]$, even though 
$\mu_{1-1/\sqrt{2}} \preccurlyeq \mu \preccurlyeq \mu_{1/\sqrt{2}}$,
as is easily verified.
We also observe that this measure does not even have full support. 
\endprocl

\procl x.ustH-ent
Let $f$ be as in \ref x.ustH/, so that $\P^f$ is the law of the horizontal
edges of a uniform spanning tree in the square lattice $\Z^2$.
It is known that the entropy of the entire uniform spanning tree measure (both
horizontal and vertical edges included) is \ref e.catalan/;
see \ref b.BurPem/. 
This can be compared to the bound on $H(\P^f)$ that
\ref p.entropy-d/ provides, together with the calculations of
\ref x.ustH-dom/, namely,
$
H(\P^f) \ge H[e^{-4\cat/\pi}] = 0.6203^+
$.
(The results are much worse if one uses the bounds of \ref b.ShiTak:II/
reported in \ref r.AMGM/ above.)
Direct calculation using cylinder events corresponding to a 4-by-4 block
gives an upper bound for the entropy $H(\P^f) \le 0.68864$.
Of course, the vertical edges of the uniform spanning tree measure have the
same entropy as the horizontal edges.
\endprocl

\procl x.ust_H-x
Let $g(x)$ be as in \ref x.ustHx/, so that
the edges of the uniform spanning tree in the plane that
lie on the $x$-axis have the law $\P^g$.
In view of \ref x.ustH_x-dom/,
\ref p.entropy-d/ implies a lower
bound of $H(\P^g) \ge H[\sqrt 2 - 1] \ge 0.67835$.
Also, direct calculation using cylinder events of length 16 gives an upper
bound $H(\P^g) \le 0.69034$.
In \ref x.ustH-x-ent/, we shall obtain more refined bounds on the entropy.
It is not hard to see that the horizontal edges associated to the $y$-axis
have law $\P^{\bfo - g}$.
\endprocl

Although not relevant for our determinantal probability measures
(where we have seen that $\strle$ and $\preccurlyeq$ are equivalent),
it is interesting to ask whether 
$\mu_p \preccurlyeq \mu \preccurlyeq \mu_{1-p}$ with $p > 0$
yields any lower bound on the entropy $H(\mu)$ for general $\mu$. 
We first ask the question whether $\mu_p \preccurlyeq \mu$
with $p > 0$ implies that
$H(\mu)> 0$ provided $\mu\neq \delta_\bfo$. The answer to
this question is affirmative
since it is known (see \ref b.furst/ for $d=1$ and \ref b.GTW/ for $d > 1$)
that zero-entropy processes are {\bf disjoint} from i.i.d.\
processes (meaning that there are no stationary couplings of them other than
independent couplings).
The next result provides an explicit lower bound on the entropy
for processes trapped between two i.i.d.\ processes.
The proof was obtained jointly with Chris Hoffman.

\procl p.hoffman
For any $d\ge 1$, if
$\mu_p \preccurlyeq \mu \preccurlyeq \mu_{1-p}$ with $p >0$, then
$$
H(\mu)\ge \max{\{a_p,b_p\}} >0
\,,
$$
where
$$
a_p:=
(1-p)\log\Big({1\over  1-p}\Big)-{1-2p \over 2}\log\Big({1\over
1-2p}\Big)
$$
and
$$
b_p:=
2(1-p)\log\Big({1\over  1-p}\Big)-(1-2p)\log\Big({1\over
1-2p}\Big)-(1-2p)\log 2.
$$
\endprocl

\procl r.log2
Observe that $b_p$ approaches $\log 2$ as $p\to 1/2$.
\endprocl

\proof
Let $X$, $Y$ and $Z$ denote processes with respective distributions 
$\mu_p,\mu$ and $\mu_{1-p}$. We construct a joining (stationary coupling)
of all three processes as follows. Consider any joining $m_1$ of $X$ and $Y$
with $X_i\le Y_i$ a.s.\ 
and any joining $m_2$ of $Y$ and $Z$ with $Y_i\le Z_i$ a.s. We now pick a 
realization for $Y$ according to $\mu$ and then choose $X$ and $Z$ 
(conditionally) independently using $m_1$ and $m_2$ respectively (this
is called the fibered product of $m_1$ and $m_2$ over $Y$). This gives
us a joining $(X',Y',Z')$ of $X$, $Y$ and $Z$.
We may now assume that $(X, Y, Z) = (X',Y', Z')$.

We now use standard facts about entropy.
First, $H(X,Z)\ge H(X,Y,Z)- H(Y)$. We next note, using
the conditional independence of $X$ and $Z$ given $Y$, that
\begineqalno
H(X,Y,Z)
&= H(Y) + H(X,Z\mid Y)=
H(Y) + H(X\mid Y)+H(Z\mid Y)
\cr&\ge  H(X)+H(Z)-H(Y)\,. 
\cr
\endeqalno
Hence
$$
H(X,Z)\ge H(X)+H(Z)-2H(Y)=2H[p]-2H(Y)\,.
$$
Since $X\le Z$, the one-dimensional marginal of $(X,Z)$ necessarily
has 3 atoms of weights $p,1-2p$ and $p$. It follows that
$$
H(X,Z)\le 2p\log\Big({1\over p}\Big)+(1-2p)\log\Big({1\over 1-2p}\Big)\,.
$$
Combining this with the previous inequality shows that 
$H(\mu)\ge a_p$. 

We modify the above proof to show that $b_p$ is also
a lower bound. Rather than using $H(X,Z)\ge H(X,Y,Z)- H(Y)$, we
use $H(X,Z)= H(X,Y,Z)- H(Y\mid X,Z)$. The earlier argument then gives
$$
H(X,Z)\ge 2H[p]-H(Y)- H(Y\mid X,Z)\,.
$$
Using the same upper bound on $H(X,Z)$ as before and the trivial upper bound
for $H(Y\mid X,Z)$ of $(1-2p)\log 2$ (obtained by noting that $Y_i$ is
determined by $X_i$ and $Z_i$ when $X_i=Z_i$) yields the lower bound of $b_p$.

The easiest way to check the strict positivity of $a_p$ is to
observe that
$$
2H[p] > 2p\log\Big({1\over p}\Big)+(1-2p)\log\Big({1\over 1-2p}\Big)
\,.
$$
This, in turn, follows from the fact that the left-hand side is the entropy of
the independent joining $\mu_p \times \mu_{1-p}$, while the right-hand side is
the entropy of the monotonic coordinatewise-independent joining of $\mu_p$
and $\mu_{1-p}$, whose existence follows from the inequality 
$p \le 1-p$.
\Qed

There are few cases where we know the entropy $H(\P^f)$ exactly.
Besides the case when $f$ is constant, we can calculate the entropy when $f$
is the reciprocal of a trigonometric polynomial of degree 1.
We shall illustrate this in the simplest case, \ref x.renewal/.
Thus, let $0 < a < 1$, $d = 1$, and $f(x) := (1 - a)^2/|e^{2\pi i x} - a|^2$.
Recall that that for any invariant measure $\mu$ in one dimension, we have
$$
H(\mu) =
\int H\big[\mu[\eta_0 = 1 \mid \eta_{-1}, \eta_{-2}, \ldots]\big]
\,d\mu(\eta_{-1}, \eta_{-2}, \ldots)
\,.
\label e.ent-past
$$
Now by \ref e.renewal-distribution/, it is easy to check that
$$
P^f[\eta_0 = 1 \mid \eta_{-1}, \eta_{-2}, \ldots] =
{(1 - a)^2 N \over N - (N - 1) a}
\,,
$$
where $N := \min \{ k\ge 1 \st \eta_{-k} = 1 \} $.
In addition, we may easily show that for $n\ge 1$,
$$
\P^f[N = n] = {(1-a)(n - (n-1)a) a^{n-1} \over 1+a}
\,.
$$
Therefore 
\begineqalno
H(\P^f) 
&=
\sum_{n=1}^\infty 
{(1-a)(n - (n-1)a) a^{n-1} \over 1+a}
\Big(
- {(1 - a)^2 n \over n - (n - 1) a} \log {(1 - a)^2 n \over n - (n - 1) a}
\cr&\hskip1in
- {n - (n - 1) a - (1 - a)^2 n \over n - (n - 1) a} 
\log {n - (n - 1) a - (1 - a)^2 n \over n - (n - 1) a}
\Big)
\cr&=
{1-a \over 1+a}
\sum_{n=1}^\infty a^{n-1}
\Big(
- (1 - a)^2 n \log [(1 - a)^2 n]
\cr&\hskip1.3in
- [a(1 - a) n + a] \log [a(1 - a) n + a]
\cr&\hskip1.3in
+ [(1-a) n + a] \log [(1-a) n + a]
\Big)
\,.
\cr
\endeqalno

In principle, one can also write an infinite series of positive terms for the
entropy of $\P^f$ when $f$ is equal to the reciprocal of any trigonometric
polynomial, since, by \ref p.regenerate/, the process $\P^f$ is then
a regenerative process.
Of course, the answer will be much more unwieldy than the above formula.
However, it can be used to get arbitrarily good lower bounds on the entropy of
any process $\P^f$, in theory.
To see this, 
we use the following lemma, which is more or less Lemma 5.10 in 
\ref b.Rudolph/. As our proof is somewhat shorter and gives a more precise
bound, we include it. See also \ref b.BurPem/, Lemma 6.2, for a similar
proof.

\procl l.rudolph
For any stationary processes $\mu$ and $\nu$, we have
$|H(\mu)-H(\nu)| \le H[\dbar(\mu,\nu)]$.
\endprocl

\proof
Consider a
stationary process $(X,Y)=\Seq{(X_i,Y_i)}_{i\in\Z^d}$,
where $X=\Seq{X_i}_{i\in\Z^d}$ has distribution $\mu$,
$Y=\Seq{Y_i}_{i\in\Z^d}$ has distribution $\nu$, and
$P(X_0\neq Y_0)=\dbar(\mu,\nu)$. 
Let $Z_i := X_i + Y_i$ mod 2. Then 
$$
H(\mu) \le H(X, Z) = H(Y, Z) \le H(\nu) + H(Z)        
\le H(\nu) + H[\dbar(\mu, \nu)]
\,.
$$
Since the same holds with $\mu$ and $\nu$ switched, the result follows.
\Qed

Thus, given any measurable $f : \T \to [0, 1]$ and any $\epsilon > 0$, define
$\delta$ to be the smallest positive number 
so that $H[\delta] = \epsilon$.
Then choose a trigonometric polynomial $T \ge 1$ such that $\int_\T |T^{-1} -
\max(f,\delta/2)| \,d\lambda_1 <\delta/2$.
We can calculate $H(\P^{1/T})$ as closely as desired.
Since $\dbar(\P^{1/T}, \P^f) < \delta$ by \ref p.L1givesdbar/,
it follows by \ref l.rudolph/ that $|H(\P^f) - H(\P^{1/T})| < \epsilon$.

Unfortunately,
this method of calculation is hopeless in practice since when $T$ has a high
degree, it will take a very long time to see a renewal, which means that
cylinder events of great length will be needed for the estimation.
When a cylinder event of length $n$ is used, one must calculate $2^n$
probabilities. Thus, this is completely impractical.

However, we now exhibit an alternative method that works extremely well in
practice, although we cannot prove that it works well {\it a priori}.
We begin with two simple examples,
$$
f(x) := \I{[0, 1/2]}(x)
\label e.deff
$$
(which also occurred in \ref x.domino/)
and 
$$
g(x) := \sin^2 \pi x = (1 - \cos 2\pi x)/2
\,.
\label e.defg
$$
Since $\widehat f(k) = \widehat g(k) = 0$ for all even $k \ne 0$,
both $\P^{f}$ and $\P^{g}$ have the property that looking at only the even
coordinates, we see independent fair coin flips, i.e., $\mu_{1/2}$.
Therefore the entropy of both processes is at least $(1/2) \log 2$.
In either case, we know of no direct method to prove that strict inequality
holds, but it does.
We first show this for $\P^g$. 
In \ref x.ustH_x-dom/ we showed that $\GM(g) = 1/4$.
By symmetry, $\GM(1-g) = 1/4$.
Since $H[1/4] = 0.56^+ > (1/2) \log 2 = 0.35^-$, we obtain by
\ref p.entropy-d/ that $H(\P^g) > (1/2) \log 2$.

In order to show that $H(\P^f) > (1/2) \log 2$ for $f$ as in \ref
e.deff/,
we need a method to
obtain more refined entropy bounds. We illustrate such a method
beginning with this simple function $g$.
While we shall explain afterwards a method to obtain results for more general
functions for the case $d=1$, this first proof contains the essential
idea of this method while at the same time relying on a more elementary 
calculation, and therefore, we feel, is worth including.

\procl p.entropyg
With $g$ as in \ref e.defg/, we have 
$$
0.63^- = {3 \over 8} H[1/4] + {5 \over 8} H[11/28]
\le H(\P^g) \le
{3 \over 8} H[7/20] + {5 \over 8} H[5/12] = 0.67^-
\,.
$$
\endprocl

\procl l.GMpert
Let $h : \T \to [0, 1]$ be a trigonometric polynomial of degree at most 1,
i.e., $\widehat h(k) = 0$ for $|k| \ge 2$.
Fix $n > 0$ and $A \subseteq \{ 1, 2, \ldots, n \}$.
For $C \subseteq \{ 1, 2, \ldots, n \}$, let $\lambda(C)$ denote the sequence
of lengths of consecutive intervals in $\{ 1, 2, \ldots, n \} \setminus C$.
Then 
$$
\P^h[\eta \restrict A \equiv 0,\, \eta \restrict (\{ 1, 2, \ldots, n \}
\setminus A) \equiv 1]
=
\sum_{C \subseteq A} (-1)^{|A \setminus C|} \prod_{i \in\lambda(C)} D_{i-1}(h)
\,.
$$
\endprocl

\proof
Whenever two sets $F_1, F_2 \subset \Z$ neither overlap nor come within
distance 1 of each other, the configurations on $F_1$ and $F_2$ are
$\P^h$-independent because $\widehat h(k) = 0$ for $|k| \ge 2$.
Therefore $\P^h[\eta \restrict (\{ 1, 2, \ldots, n \} \setminus C) \equiv 1] =
\prod_{i \in\lambda(C)} D_{i-1}(h)$. Now the desired formula follows from the
inclusion-exclusion principle.
\Qed

\proofof p.entropyg
Write $d_n := D_{n-1}(g)$.
Recall that 
$$
H(\P^g) =
\int H\big[\P^g[\eta_0 = 1 \mid \eta_{-1}, \eta_{-2}, \ldots]\big]
\,d\P^g(\eta_{-1}, \eta_{-2}, \ldots)
\,.
\label e.ent-past
$$

By the negative association property of $\P^g[\;\cdot \mid \eta_{-1} =
\eta_{-2} = 1]$, we have that on the event $ \{ \eta_{-1} = \eta_{-2} = 1 \}$,
\begineqalno
\P^g[\eta_0 = 1 \mid \eta_{-1}, \eta_{-2}, \ldots]
&\ge
\lim_{n \to\infty} \P^g[\eta_0 = 1 \mid \eta_{-1} = \eta_{-2}= \cdots = \eta_{-n} = 1]
\cr&=
\lim_{n \to\infty} {\P^g[\eta_0 = \eta_{-1} = \eta_{-2}= \cdots = \eta_{-n} = 1] \over
\P^g[\eta_{-1} = \eta_{-2}= \cdots = \eta_{-n} = 1]}
\cr&=
\lim_{n \to\infty} d_{n+1}/d_n
=
\GM(g)
=
1/4
\,.
\endeqalno
Likewise, on the same event, we have 
\begineqalno
\P^g[\eta_0 = 1 \mid \eta_{-1}, \eta_{-2}, \ldots]
&\le
\lim_{n \to\infty} \P^g[\eta_0 = 1 \mid \eta_{-1} = \eta_{-2}= 1,\, \eta_{-3} = 
\cdots = \eta_{-n} = 0]
\cr&=
\lim_{n \to\infty} \P^g[\eta_0 = 0 \mid \eta_{-1} = \eta_{-2}= 0,\, \eta_{-3} = 
\cdots = \eta_{-n} = 1]
\cr&\hskip1in\hbox{[by symmetry]}
\cr&=
1 - \lim_{n \to\infty} \P^g[\eta_0 = 1 \mid \eta_{-1} = \eta_{-2}= 0,\, \eta_{-3} = 
\cdots = \eta_{-n} = 1]
\cr&=
1 - \lim_{n \to\infty} {\P^g[\eta_{-1} = \eta_{-2}= 0,\, \eta_0 = \eta_{-3} = 
\cdots = \eta_{-n} = 1] \over
\P^g[\eta_{-1} = \eta_{-2}= 0,\, \eta_{-3} = \cdots = \eta_{-n} = 1]}
\cr&=
1 - \lim_{n \to\infty} 
{d_1 d_{n-2} - d_2 d_{n-2} - d_1 d_{n-1} + d_{n+1}
\over
d_{n-2} - d_1 d_{n-2} - d_{n-1} + d_{n}}
\cr&\hskip1in\hbox{[by \ref l.GMpert/]}
\cr&=
1 - 
{d_1 - d_2 - d_1 \GM(g) + \GM(g)^3
\over
1 - d_1 - \GM(g) + \GM(g)^2}
\cr&=
7/20
\endeqalno
since $d_1 = 1/2$ and $d_2 = 3/16$.
We may conclude that on the event $ \{ \eta_{-1} = \eta_{-2} = 1 \}$, we have 
$$
H\big[\P^g[\eta_0 = 1 \mid \eta_{-1}, \eta_{-2}, \ldots]\big]
 \in \big[H[1/4], H[7/20]\big]
\,.
$$
By symmetry, the same holds on the event $ \{ \eta_{-1} = \eta_{-2} = 0 \}$.

Similarly, we have that on the event $ \{ \eta_{-1} = 1,\, \eta_{-2} = 0 \}$,
\begineqalno
\P^g[\eta_0 = 1 \mid \eta_{-1}, \eta_{-2}, \ldots]
&\ge
\lim_{n \to\infty} \P^g[\eta_0 = 1 \mid \eta_{-2} = 0,\,
\eta_{-1} = \eta_{-3}= \cdots = \eta_{-n} = 1]
\cr&=
\lim_{n \to\infty} {\P^g[\eta_{-2} = 0,\,
\eta_0 = \eta_{-1} = \eta_{-3}= \cdots = \eta_{-n} = 1] \over
\P^g[\eta_{-2} = 0,\, \eta_{-1} = \eta_{-3}= \cdots = \eta_{-n} = 1]}
\cr&=
\lim_{n \to\infty}
{d_2 d_{n-2} - d_{n+1}
\over
d_1 d_{n-2} - d_n}
\cr&=
{d_2 - \GM(g)^3
\over
d_1 - \GM(g)^2}
\cr&=
11/28
\,.
\endeqalno
Likewise, on the same event, we have 
\begineqalno
\P^g[\eta_0 = 1 \mid \eta_{-1}, \eta_{-2}, \ldots]
&\le
\lim_{n \to\infty} \P^g[\eta_0 = 1 \mid \eta_{-1} = 1,\, \eta_{-2}= \eta_{-3} = 
\cdots = \eta_{-n} = 0]
\cr&=
\lim_{n \to\infty} \P^g[\eta_0 = 0 \mid \eta_{-1} = 0,\, \eta_{-2}= \eta_{-3} = 
\cdots = \eta_{-n} = 1]
\cr&=
1 - \lim_{n \to\infty} \P^g[\eta_0 = 1 \mid \eta_{-1} = 0,\, \eta_{-2}= \eta_{-3} = 
\cdots = \eta_{-n} = 1]
\cr&=
1 - \lim_{n \to\infty} {\P^g[\eta_{-1} = 0,\, \eta_0 = \eta_{-2} = 
\cdots = \eta_{-n} = 1] \over
\P^g[\eta_{-1} = 0,\, \eta_{-2} = \cdots = \eta_{-n} = 1]}
\cr&=
1 - \lim_{n \to\infty} {d_1 d_{n-1} - d_{n+1} \over d_{n-1} - d_n}
\cr&=
1 - {d_1 - \GM(g)^2 \over 1 - \GM(g)}
\cr&=
5/12
\,.
\endeqalno
We may conclude that on the event $ \{ \eta_{-1} = 1,\, \eta_{-2} = 0 \}$, we have 
$$
H\big[\P^g[\eta_0 = 1 \mid \eta_{-1}, \eta_{-2}, \ldots]\big]
 \in \big[H[11/28], H[5/12]\big]
\,.
$$
By symmetry, the same holds on the event $ \{ \eta_{-1} = 0,\, \eta_{-2} =
1 \}$.

Putting all these bounds into \ref e.ent-past/ gives the claimed bounds on the
entropy of $\P^g$.
\Qed

For more general functions, we can get lower bounds on the entropy by using
the asymptotics of \ref b.BumpDiaconis/ or of \ref b.TW:like/,
which serve as replacements for our use of the special \ref l.GMpert/.
In fact, we use an extension of the formula of \ref b.TW:like/ due to
\ref b.Lyons:szego/.
Denote by $\nu_{f}$ the measure on $2^\N$ obtained by the
limiting condition $\lim_{n \to\infty} \P^f[\;\cdot \mid \eta_{-1} = \eta_{-2}
= \cdots = \eta_{-n} = 1]$, which exists by \ref p.JNRD/.
If $\GM(f) > 0$, define
$$
\Phi_f(z) := \exp {1 \over 2} \int_{\T} {e^{2\pi i t} + z
\over e^{2\pi i t} - z} \log f(t) \,d\lambda_1(t)
\label e.defPhi_f
$$
for $|z| < 1$.
The {\bf outer} function
$$
\varphi_f(t) := \lim_{r \uparrow 1} \Phi_f(r e^{2 \pi i t})
\label e.defvarphi_f
$$
exists for $\lambda_1$-a.e.\ $t \in \T$ and satisfies $|\varphi_f|^2 = f$
$\lambda_1$-a.e. 
This limit also holds in $L^1(\T)$.
See \ref b.Rudin:RCA/, Theorem 17.11, p.~340, and Theorem 17.16, p.~343.
As is easily verified,
$$
\log \Phi_f(z) = \widehat F(0) + 2 \sum_{k=1}^\infty \widehat F(k) z^k
\label e.preexp
$$
when $F := (1/2) \log f$.
The mean value theorem for analytic functions and the $L^1(\T)$
convergence above show that
$$
\widehat {\varphi_f}(0) = \Phi_f(0) = \sqrt{\GM(f)}
\,.
$$
Define $\Phi_f := \varphi_f := \bfz$ if $\GM(f) = 0$.
An analytic trigonometric polynomial is (a constant times) an outer
function iff its extension to the unit disc has no zeroes in the open disc.
[On the one hand, the extension of 
no outer function has any zeroes in the open disc since it
is an exponential.  On the other hand, by factoring a polynomial, one
has to check merely that $z \mapsto z-\zeta$ is a constant times
an outer function for
$|\zeta| \ge 1$. Indeed, $\Phi_f(z) = (-|\zeta|/\zeta)(z -\zeta)$
for $f(t) := |e^{2\pi i t} -\zeta|^2$, as can be seen by the Poisson integral
formula for $|\zeta|>1$ (\ref b.Rudin:RCA/, Theorem 11.9, p.~235) and by a
limiting procedure for $|\zeta| = 1$.]
The results of \ref b.Lyons:szego/ show that \ref t.condonB/ reduces to
the following formula:

\procl t.szego
For any measurable $f : \T \to [0, 1]$, the measure $\nu_{f}$ is equal to
the determinantal probability measure corresponding to the positive
contraction on $\ell^2(\N)$ whose $(j, k)$-matrix entry is 
$$
\sum_{l=0}^{j \wedge k} \overline{\widehat{\varphi_f}(j-l)}
\widehat{\varphi_f}(k-l)
\,.
$$
\endprocl

This can be substituted in the appropriate places in the proof of \ref
p.entropyg/.
For example, on the event $ \{ \eta_{-1} = \eta_{-2} = 0 \}$, one has that 
\begineqalno
\P^f[\eta_0 = 1 \mid \eta_{-1}, \eta_{-2}, \ldots]
&\ge
\lim_{n \to\infty}
\P^f[\eta_0 = 1 \mid \eta_{-1} = \eta_{-2}= 0,\, \eta_{-3} = \eta_{-4} =
\cdots = \eta_{-n} = 1]
\cr&=
\nu_{f}[\eta_2=1, \eta_1=\eta_0=0]/ \nu_{f}[\eta_1=\eta_0=0]
\,.
\endeqalno
One can use $\nu_{\bfo-f}$ in a similar fashion for estimating such
probabilities above by $\eta$ that terminate in repeating 0s, rather than in
repeating 1s.
For example, if 
$$
f := |1+\ee_{1} + \ee_{2}|^2/9
\,,
$$
then $\varphi_f = (1+\ee_{1} + \ee_{2})/3$ since the polynomial
$(1+z+z^2)/3$ has no zeroes in the open unit disc.
Similarly, one can show that 
$$
\varphi_{\bfo-f} = {\sqrt 6 + \sqrt 2 \over 6} - {\sqrt 2 \over 3} \ee_{1} 
+ {\sqrt 2 - \sqrt 6 \over 6} \ee_{2}
\,.
$$
With these and \ref t.szego/, one can show by following the method of proof of
\ref p.entropyg/ that 
$$
0.53 \le H(\P^f) \le 0.61
\,.
$$
This method of estimation is relatively easy to program on a computer;
decomposing by cylinder events of length 8 instead of the length 2 used in the
proof of \ref p.entropyg/ and using Mathematica, we find that 
$$
0.601992  \le H(\P^f) \le 0.602433
\,.
$$
Similarly, using cylinder events of length 15, we obtain that 
$$
0.65907716
\le H(\P^g) \le
  0.65907733
\label e.entropy-g
$$
for the function $g$ in \ref e.defg/.

We now return to the problem of showing that
$H(\P^f) > (1/2) \log 2$ for $f$ as in \ref e.deff/.

\procl p.entrindic
If $f := \I{[0, 1/2]}$, then $H(\P^f) > (1/2) \log 2$.
\endprocl

\proof
Let $\widetilde{f}:= .99\I{[0, 1/2]} +.01\I{(1/2,1)}$.
\ref p.L1givesdbar/ tells us that
$\dbar(\P^f, \P^{\widetilde{f}}) \le .01$, and hence from
\ref l.rudolph/,
$H(\P^f) \ge H(\P^{\widetilde{f}})- H[.01]$. 
Now $H[.01] < 0.0561$.
Write $\delta := .01$ and
$L := \log (\delta ^{-1} - 1) = \log 99$.
A simple integration shows that 
$$
\skew6\widehat {\widetilde f}(k) =
\cases{(1/2) &if $k = 0$,\cr
\displaystyle{(-1)^k - 1 \over 2k\pi} i (1 - 2\delta)
&if $k \ne 0$\cr}
$$
(where $i := \sqrt{-1}$).
To find $\varphi_{\widetilde f}$, we proceed as follows.
Let $F := (1/2) \log \widetilde f$.
Then
$$
\widehat F(k) =
\cases{(1/4) \log \delta(1-\delta)  &if $k = 0$,\cr
\displaystyle{(-1)^k - 1 \over 4k\pi} i L
&if $k \ne 0$.\cr}
$$
By using \ref e.preexp/ and exponentiating, we obtain
$$
\Phi_{\widetilde f}(z)
=
[\delta(1-\delta)]^{1/4}
\left(
1
-
{i \over \pi} L z
-
{1 \over 2\pi^2} L^2 z^2
-
{i \over 6\pi^3} (2\pi^2 - L^2) L z^3
+
\cdots
\right)
\,.
$$
Since the coefficients of this power series are the Fourier coefficients of
$\varphi_{\widetilde f}$,
we have a procedure to compute the Fourier coefficients
$\varphi_{\widetilde f}$.
Since $\widetilde f(1-x) = 1-\widetilde f(x)$, we have
$\overline{\Phi_{\widetilde f}(\overline z)} =
\Phi_{\bfo-\widetilde f}(z)$. 
Thus, the $k$th Fourier coefficient of $\varphi_{\bfo - \widetilde f}$ equals
the complex conjugate of the $k$th Fourier coefficient of $\varphi_{\widetilde
f}$.
Using the method of \ref p.entropyg/ with cylinder events of length 3,
which requires knowledge of Fourier coefficients only for $|k| \le 2$, we
obtain that $H(\P^{\widetilde f}) \ge 0.4105$.
Therefore $H(\P^f) > 0.4105 - 0.0561 = 0.3544 > (1/2) \log 2$.
(If we use cylinder events of length 12, then we obtain the bound $H(\P^f) >
0.4442$. However, we believe that the true entropy is significantly larger
still.)
\Qed

\procl x.depends
We can also use these bounds to prove that $H(\P^f)$ does {\it not\/}
depend {\it only\/} on the distribution of $f$.
For example, consider the function 
$$
f(t) := \sin^2(\pi t/2)
$$
on $[0, 1]$, which has the same distribution as the function $g$ of \ref
e.defg/. 
An elementary integration shows that 
$$
\widehat f(k) = \cases{1/2 &if $k=0$,\cr
\displaystyle {2k i \over (2k-1)(2k+1) \pi} &if $k \ne 0$\cr}
$$
(where $i := \sqrt{-1}$).
To find $\varphi_f$, we proceed as follows.
Write $F := (1/2) \log f$.
The real parts of the integrals giving $\widehat F(k)$ can be found in
standard tables, while the imaginary parts are derived in \ref b.LPR/,
giving
$$
\widehat F(k) =
\cases{-\log 2  &if $k=0$,\cr
       \displaystyle -{1 \over 4k} + {i \over k \pi} \sum_{j=1}^k {1 \over 2j-1}
                &if $k > 0$.\cr}
$$
As in the proof of \ref p.entrindic/, we compute
\begineqalno
\Phi_f(z) =
1/2 &+( -1/4 + i/\pi) z +(-1/16 + i/6\pi - 1/\pi^2) z^2
\cr&  -(1/32 - 19 i/360 \pi + 5/6 \pi^2 + 2 i/3\pi^3) z^3
\cr&  +(-5/256 + 89 i/5040\pi - 27/40 \pi^2 - i/\pi^3 + 1/3 \pi^4) z^4
  + \cdots
\,,
\endeqalno
whose coefficients are the Fourier coefficients of $\varphi_f$.
Since $1-\sin^2(\pi t/2) = \sin^2(\pi (1-t)/2)$, it follows from the
definition \ref e.defPhi_f/ that $\overline{\Phi_{f}(\overline z)} =
\Phi_{\bfo-f}(z)$. Since the power series coefficients of $\Phi_{\bfo - f}$
are the Fourier coefficients of $\varphi_{\bfo - f}$, we obtain that the $k$th
Fourier coefficient of $\varphi_{\bfo - f}$ equals the complex conjugate of
the $k$th Fourier coefficient of $\varphi_f$.
Using the method of \ref p.entropyg/ with cylinder events of length 8,
which requires knowledge of Fourier coefficients only for $|k| \le 7$, we
obtain that 
$$
0.659648
\le H(\P^f) \le
0.684021
\,.
$$
Comparing with \ref e.entropy-g/, we see that $H(\P^f) > H(\P^g)$.
\endprocl

Note that even when neither $\widehat f$ nor $\widehat{\varphi_f}$ can be
found explicitly, they can always be found by numerical integration and
exponentiation, and then one can follow the procedure we have used in 
this last example. Indeed, that will be done for part of our final example.

\procl x.ustH-x-ent
As in \ref x.ust_H-x/, let 
$$
g(x) := 
{\sin \pi x \over \sqrt{1+\sin^2 \pi x}}
\,.
$$
Then the edges of the uniform spanning tree measure in the plane that lie on
the $x$-axis have the law $\P^g$.
Write $g_1(x) := \sin^2 \pi x$.
The calculations in \ref x.ustHx/ show that $\Phi_{g_1}(z) = (1-z)/2$ and
$\Phi_{\bfo+g_1}(z) = [\sqrt 2 + 1 - (\sqrt 2 - 1)z]/2$, whence
$$
\Phi_g(z)
=
\sqrt{\Phi_{g^2}(z)}
=
\sqrt{{\Phi_{g_1} \over \Phi_{\bfo+g_1}}}
=
\sqrt{1-z \over \sqrt 2 + 1 - (\sqrt 2 - 1)z}
\,.
$$
Expansion of $\Phi_g(z)$ in a Maclaurin series gives 
$$
\frac{1}{({1 + {\sqrt{2}}})^{1/2}} - \frac{z}{( 1 + {\sqrt{2}} )
^{\nonfrac{3}{2}}} - \frac{(-1+2\sqrt{2})\,z^2}{2\,(1 +\sqrt{2})^{5/2}} -
\frac{( 5 - 2\,{\sqrt{2}} ) \,z^3}{2\,( 1 + {\sqrt{2}}
) ^{\nonfrac{7}{2}}} - \frac{( -27 + 28\,{\sqrt{2}} )
\,z^4}{8\,( 1 + {\sqrt{2}} ) ^{\nonfrac{9}{2}}} +
\cdots
\,,
$$
which tells us $\widehat{\varphi_g}$.
The transfer currents, which can be calculated by the method in \ref
b.BurPem/, tell us $\widehat g$; for example, $\widehat g(k)$ for $k=0, 1, 2,
3, 4$ is 
$$
\frac{1}{2}\,,\qquad \frac{1}{2} - \frac{2}{\pi
}\,,\qquad\frac{5}{2} - \frac{8}{\pi }\,,\qquad\frac{25}{2} - \frac{118}{3\,\pi
}\,,\qquad
\frac{129}{2} - \frac{608}{3\,\pi }
\,.
$$
We use numerical integration to find $\widehat{\varphi_{\bfo-g}}$.
Then if we use cylinder events of length 8, we find that 
$$
0.69005 \le H(\P^g) \le 0.69013
\,.
$$
It is interesting how close this is to $\log 2 = 0.69315^-$.
\endprocl

We close our treatment of entropy with some elementary observations.
If $f_n \to 1/2$ weak${}^*$ (meaning that we have weak${}^*$
convergence of the measures having these densities),
must $H(\P^{f_n}) \to \log 2$?
The answer is ``no", as we now demonstrate.
Given any $f : \T \to [0, 1]$ and any integer $n$, let $\mult{n} f$ denote the
function 
$$
\mult{n} f(x) := f(nx)
\,.
$$
Let $\Seq{\eta^{(n)}_k \st k \in \Z}$ have the distribution $\P^{\mult{n}f}$.
Since $f(x) = \sum_{k \in \Z} \widehat{f}(k) e^{2\pi i k x}$ in $L^2(\T)$, the
Fourier expansion of $\mult{n}f$ is $\mult{n}f(x) = \sum_{k \in \Z}
\widehat{f}(k) e^{2\pi i k n x}$ for $n \ne 0$.
In particular, for any $n \ne 0$ and any $r$, the processes
$\Seq{\eta^{(n)}_{nk+r} \st k \in \Z}$ 
each have distribution $P^f$ and are independent of each other 
as $r$ ranges from 0 to $n-1$.
Therefore $H(\P^{\mult{n}f}) = H(\P^f)$.
Note that unless $f$ is constant, 
$H(\P^{\mult{n}f}) = H(\P^f) < H\big[\widehat f(0)\big]$, even
though $\mult{n}f$ tends weak$^*$ to the constant function $\widehat f(0)
\cdot \bfo$.
One can show that a similar phenomenon holds in higher dimensions, that is,
if $f : \Td \to [0, 1]$ and $A: \Td \to \Td$ is a group epimorphism, then
$H(\P^f) = H(\P^{f \circ A})$.

\bsection {Phase Multiplicity}{s.sk}

In this section, we classify exactly the set of functions $f$ for 
which $\P^f$ satisfies a
strong full $K$ property. This property, which we now describe, is an 
essential strengthening of the 
usual Kolmogorov or $K$ property. One of the reasons this property 
is interesting is that it is closely connected to the notion
of multiplicity for Gibbs states in 
statistical mechanics; in particular, it corresponds to 
{\it uniqueness}. In the next section, we
classify exactly the set of functions $f$ for which $\P^f$ satisfies a
(1-sided) strong $K$ property.

{\bf In this section and the next, we always assume that $f$ is not}
{\bf identically  0 nor identically 1.}

To begin, we first recall the $K$ property for stationary processes indexed 
by $\Z$. 
There are many equivalent formulations, of which we choose an appropriate one.

For $S \subseteq \Z^d$, write $\F(S)$ for the $\sigma$-field on $2^{\Z^d}$
generated by $\eta(e)$ for $e \in S$.

\procl d.K
A translation-invariant probability measure $\mu$ on
$2^{\Z}$ is {\bf $K$\/} (or a {\bf Kolmogorov automorphism})
if for any finite cylinder event $E$ and 
for all $\epsilon>0$ there exists an $N$ such that 
$$
\mu\Big[ \big|\mu\big(E \mid \F(\OC{-\infty,-N}\cap \Z)\big)-\mu(E)\big| \ge \epsilon \Big] \le
\epsilon
\,.
$$
\endprocl

It is well known that the above definition of $K$ is equivalent 
to having a trivial (1-sided) tail $\sigma$-algebra in the sense that
the $\sigma$-algebra $\bigcap_{m \ge 1} \F\big(\OC{- \infty,
-m}\big)$ is trivial
(see page 120 of \ref b.Georgii/ for a version of this).
The $K$ property is known to be an isomorphism invariant and is
also known to be equivalent to the
property that all nontrivial factors have strictly positive entropy.
The latter implies that a process is $K$ iff its time reversal is $K$,
something which is not immediate from the definition.

In order to give a complete discussion of the points that we wish to make,
we need to introduce three further properties, which are respectively
called the full $K$ property, the strong full $K$ property, and the
(1-sided) strong $K$ property, the last property given only for $d=1$.

While the notion of the $K$ property generalizes 
to $\Z^d$ (see \ref b.Conze/),
there is a slight strengthening of the definition that 
has a more aesthetic extension to
$\Z^d$. To give this, define $B^d_n:=[-n,n]^d\cap\Z^d$.

\procl d.2K
A translation-invariant probability measure $\mu$ on
$2^{\Z^d}$ is {\bf full $K$\/} if for any finite cylinder event $E$ and 
for all $\epsilon>0$ there exists an $N$ such that 
$$
\mu\Big[ \big|\mu\big(E \mid \F((B^d_N)^c)\big)-\mu(E)\big| \ge \epsilon \Big] \le
\epsilon
\,.
$$
\endprocl

Analogously to an earlier statement, the
full $K$ property is equivalent (see again page 120 of \ref b.Georgii/)
to having a 
{\bf trivial full tail}, which means that 
$\bigcap_{m\ge 1} \F\big((B^d_m)^c\big)$ is trivial.
It is well known that $K$ does not imply 
full $K$ even for $d=1$ (see, for example, the bilaterally deterministic
Bernoulli shift processes constructed in \ref b.OW:bilat/).
However, for Markov random fields, the two notions are equivalent
(see \refbmulti{HS,HS2}).
\ref b.L:det/ proved that all (not necessarily $\Z^d$-invariant) 
determinantal probability measures satisfy
the full $K$ property. 
For Gibbs states arising in statistical mechanics, this property 
is equivalent
to an ``extremality'' property of the Gibbs state 
(see page 118 of \ref b.Georgii/ for precisely 
what this means and this equivalence). 

The following definitions strengthen the $K$ property in an even more 
essential way. We begin with the full version, which seems to us more 
natural.

\procl d.2sK
A translation-invariant probability measure $\mu$ on
$2^{\Z^d}$
is {\bf strong full $K$\/} if for any finite cylinder event $E$ and 
for all $\epsilon>0$ there exists an $N$ such that 
$$
\big|\mu\big(E \mid \F((B^d_N)^c)\big)-\mu(E)\big| <\epsilon
\quad \mu\hbox{-a.e.}
$$
\endprocl

For the full $K$ property, ``$\mu$-most'' conditionings far away have
little effect on a ``local event'', while in the strong full $K$ case,
{\it all\/} conditionings far away have little effect on a ``local event''.
This is a substantial difference. An example that illustrates this 
difference is the Ising model in $\Z^2$. 
The plus state for the Ising model at high temperatures
is strong full $K$, while at low temperatures,
it is full $K$ (because it is extremal) but not strong full $K$.

Finally, if $d=1$, we define
strong $K$ in the following way (which the reader can presumably
anticipate). One extension of this definition 
to $\Z^d$ (among various possible) will be given in \ref s.1side/.

\procl d.1sK
A translation-invariant probability measure $\mu$ on $2^{\Z}$ 
is {\bf strong $K$\/} if for any finite cylinder event $E$ and 
for all $\epsilon>0$ there exists an $N$ such that 
$$
\big|\mu\big(E \mid \F(\OC{-\infty,-N}\cap \Z)\big)-\mu(E)\big| <\epsilon
\quad\mu\hbox{-a.e.}
$$
\endprocl

We note that this definition is closely related to, but weaker than,
the $\psi$-mixing property (see \ref b.Bradley/). The difference is that the
event $E$ here is specified in advance, rather than having a uniformity 
for all events $E$, provided only that they depend 
on the random variables
with positive index. If we extend the notion of strong $K$ 
in the obvious way to the case when $\{0,1\}$ is replaced by a countable set,
then the example in \ref b.Bradley/ of a process that is 
$\psi$-mixing, but whose time reversal is not $\psi$-mixing, yields an example
of a strong $K$ process whose time reversal is not.
This cannot occur for our measures $\P^f$ since they are all time reversible.

The following development has an analogy with the ``plus and minus states''
in the Ising model from statistical mechanics; however, such knowledge 
is not needed by the reader.

Consider a function $f:\Td\to [0,1]$ and consider the corresponding probability measure
$\P^f$ on $2^{\Z^d}$. We shall define a probability measure
$(\P^f)^+$ on $2^{\Z^d}$ which will be ``$\P^f$ conditioned 
on all 1s at $\infty$'';
more specifically (but still not precisely), we want to define $(\P^f)^+$ by
$$
\Pfp := \lim_{n\to \infty} \P^f[\;\cdot \mid  \eta\equiv 1 \hbox{ on }
(B^d_n)^c]\,.
$$
To make sure this is well defined, we proceed in stages. 
Let
$$
(\Pfp)_n := \lim_{k\to\infty} \P^f[\;\cdot \mid  \eta\equiv 1 \hbox{ on
}B^d_{n+k}\backslash B^d_n]\,.
$$
This limit is taken in the weak$^{*}$-topology. \ref p.JNRD/ implies
that, when restricted to $B^d_n$, 
this sequence is stochastically decreasing and hence necessarily converges. 
$(\Pfp)_n$ is therefore well defined. One next defines
$$
\Pfp := \lim_{n\to \infty} (\Pfp)_n
\,,
$$
where again the limit is taken in the weak$^{*}$-topology. 
\ref p.JNRD/ again implies
that for fixed $k$,
$(\Pfp)_n$ restricted to $B^d_k$, is, for $n>k$,
stochastically increasing in $n$ and hence converges.
This implies that its limit $\Pfp$ is well
defined and completes the definition of $\Pfp$.
The stochastic monotonicity results also imply that 
$\Pfp \preccurlyeq  \P^f$ and that
$\Pfp$ is translation invariant.

In exactly the same way (using 0 instead of 1 boundary conditions), one defines
$\Pfm$, which satisfies $\P^f  \preccurlyeq \Pfm$. Analogy with the
ferromagnetic Ising model in statistical mechanics leads us to the
following definition.

\procl d.pt
The probability measure $\P^f$ has {\bf phase multiplicity} if $\Pfm\neq \Pfp$, and otherwise has {\bf phase uniqueness}.
\endprocl

\procl l.pt
$\P^f$ has phase uniqueness if and only if it is strong full $K$.
\endprocl

The reason this is true is that \ref p.JNRD/ implies
that the most
extreme boundary conditions are when all 1s or all 0s are used and the
measures corresponding to any other boundary conditions are
``stochastically trapped'' between these two special cases. 
A detailed proof follows straightforwardly from the stochastic 
monotonicity arguments above and
is left to the reader.

We shall now
obtain necessary and sufficient conditions for the equality of $\P^f$ and
$\Pfp$. From these, a necessary and sufficient condition for the strong full $K$ property will easily emerge.

Let $\calT$ denote the set of trigonometric polynomials on $\T^d$.
Let $L^2(1/f)$ denote the set 
$$
\Big\{h:\Td\to \C \st \int_\Td {|h|^2 \over f} \, d\lambda_d <\infty\Big\}\,. 
$$
Here we use the convention that $0/0:=0$. Note also that
$h$ needs to vanish where $f$ does.

Our main result on phase multiplicity is the following.

\procl t.sn
Assume that $f:\Td\to [0,1]$ is measurable. 
The following are equivalent. 
\beginitems
\itemrm{(i)} $\Pfp=\P^f$; 
\itemrm{(ii)} $f$ is in the closure in $L^2(1/f)$ of $\calT\cap L^2(1/f)$;
\itemrm{(iii)} There exists a nonzero trigonometric polynomial $T$ such that 
${|T|^2 \over f}\in L^1(\Td)$; i.e., $\calT\cap L^2(1/f)\neq {0}$.
\enditems
Moreover, if $\Pfp\neq \P^f$, then $\Pfp= \delta_\bfz$.
\endprocl

\proof
We shall first show that (i) and (ii) are equivalent and then
that (ii) and (iii) are equivalent.

{\bf (i) implies (ii):} Let $u_n$ be the element in $[B]_f$ achieving the
minimum in \ref e.closest/ for $B := (B^d_n)^c$.
Then 
$$
\|\bfo\|^2_f = \|\bfo - u_n\|^2_f + \|u_n\|^2_f
\,.
$$
By (i) and \ref e.closest/, we have $\|\bfo - u_n\|_f \to \|\bfo\|_f$ as $n
\to\infty$, or in other words, $\|u_n\|_f \to 0$ as $n \to\infty$.
Furthermore, $\bfo - u_n \perp [(B^d_n)^c]_f$ in $L^2(f)$, which is the same
as $(\bfo - u_n) f$ being a trigonometric polynomial $T_n$ with $\widehat
{T_n}$ supported in $B^d_n$.
We have 
$$
\|f - T_n\|^2_{(1/f)}
=
\int |f - T_n|^2 {1 \over f} \,d\lambda_d
=
\int \left|\bfo - {T_n \over f}\right|^2 f \,d\lambda_d
=
\int |u_n|^2 f \,d\lambda_d
=
\|u_n\|^2_f
\label e.twoways
$$
tends to 0 as $n \to\infty$, which proves (ii).

{\bf (ii) implies (i):}
Given $\epsilon > 0$, let
$T$ be a trigonometric polynomial with $\|f - T\|_{(1/f)} < \epsilon$.
Let $n$ be such that $\widehat T$ is supported in $B^d_n$.
Define $v := T/f$, which is in $L^2(f)$.
Then $v \perp [(B^d_n)^c]_f$ in $L^2(f)$ and $\|\bfo - v\|_f = \|f -
T\|_{(1/f)} < \epsilon$, so 
$\|P_{[(B^d_n)^c]_f}^\perp \bfo\|_f^2 > {\|\bfo\|^2_f -\epsilon^2}$.
Combining this with \ref e.byproj/, we see that
$\Pfp[\eta(\bfz) = 1] > \P^f[\eta(\bfz) = 1] - \epsilon^2$. Since
this holds for any $\epsilon > 0$, we obtain that $\Pfp[\eta(\bfz) = 1] \ge
\P^f[\eta(\bfz) = 1]$. However, 
since $\Pfp \preccurlyeq \P^f$ and both are translation invariant, 
the two measures are actually equal by \ref l.weak/, which proves (i).

{\bf (ii) implies (iii):} 
This is immediate.

{\bf (iii) implies (ii):}
Assume there exists a nonzero
trigonometric polynomial $T \in L^2(1/f)$. 
Let $\epsilon > 0$. 
Choose $B\subseteq \Td$ such that
$$
\int_B \left({|T|^2 \over f} + 2|T| +f\right) < \epsilon 
$$
and $B$ is an open set containing the zero set of $T$ (which is a set of
measure 0, as remarked in the proof of \ref t.support/). Let 
$g := T  \I B + f  \I {B^c}$. Then $g/T \in L^\infty(\Td)$ since $T$ 
is bounded away from 0 on $B^c$ and $0\le f \le 1$. We can now choose 
trigonometric 
polynomials $p_n$ on $\Td$ such that $p_n \to g/T$ a.e.\ 
on $\Td$ with $\|p_n\|_\infty \le \|g/T\|_\infty$ for all $n$ (just use
$p_n:=(g/T) * K^d_n$, where $K^d_n$ is, as before, the Fej\'er kernel for
$\Td$).  
Since $|T|^2/f \in L^1(\Td)$, it follows from the 
Lebesgue dominated convergence theorem that
$$
\lim_{n\to\infty}\int_\Td {|T|^2 \over f}  \left|{g \over T}-p_n\right|^2 =0\,.
$$
Hence there exists a trigonometric polynomial $h$ such that
$$
\int_\Td {|T|^2 \over f}  \left|{g \over T}-h\right|^2 < \epsilon\,.
\label e.1
$$
Minkowski's inequality (applied to $L^2(|T|^2/f)$) yields
$$
\left(\int_\Td {|T|^2 \over f}  \left|{f \over T}-h\right|^2\right)^{{1\over
2}} \le
\left(\int_\Td {|T|^2 \over f}  \left|{f \over T}-{g \over
T}\right|^2\right)^{{1\over 2}} +
\left(\int_\Td {|T|^2 \over f}  \left|{g \over T}-h\right|^2\right)^{{1\over
2}}\,.
$$
The second summand is at most $\epsilon^{{1/ 2}}$ by \ref e.1/. 
The first summand is 
$$
\left(\int_B {|T|^2 \over f}  \left|{f \over T}-1\right|^2\right)^{{1\over 2}}
\le 
\left(\int_B f+ 2|T| + {|T|^2 \over f}\right)^{{1\over 2}} <\epsilon^{{1\over
2}}
$$
by choice of $B$.
Hence 
$\int_\Td |f-Th|^2 {1 \over f} \, dx < 4\epsilon$. Since $Th$ is a 
trigonometric polynomial, (ii) is proved.

Finally, assume that $\Pfp\neq \P^f$. Since (iii) fails,
$[(B^d_n)^c]_f^\perp = \bfz$ for all $n$, which means by \ref e.byproj/
that $\Pfp[\eta(\bfz) = 1] = 0$, or, in other words, $\Pfp = \delta_\bfz$.
\Qed

\procl x.continuous-finite
If $f:\T\to [0,1]$ is continuous and has a finite number of 0s with $f$
approaching each of these 0s at most polynomially quickly, then $\Pfp=\P^f$
since it is easy to construct a trigonometric polynomial $T$ with
the same 0s as $f$ and approaching 0 at least as quickly as $f$
(for example, $T$ could be of the form 
$\prod_{i=1}^k \sin^n 2\pi(x-a_i)$).
In particular, if $f$ is real analytic and not $\bfz$, 
then $\Pfp = \P^f$. 
\endprocl

\procl x.trivial
If $f$ vanishes on a set of positive measure, then $\Pfp=\delta_\bfz$.
\endprocl

\procl x.cont
If $f:\T\to [0,1]$ is a continuous function with a single zero at $x_0\in \T$
and $f(x) = e^{-1/|x-x_0|}$ in some neighborhood of $x_0$,
then $\Pfp=\delta_\bfz$. Indeed, there is no nonzero 
trigonometric polynomial $T$ with
$|T|^2/f\in L^1(\T)$ since the rate at which a trigonometric polynomial 
approaches 0 is at most polynomially quickly.  
\endprocl

\procl x.depends-distribution
We give an example to show that the property $\Pfp = \P^f$ does {\it not\/}
depend only on the distribution of $f$, even among real-analytic functions
$f$.
In 2 dimensions, the function 
$$
f(x, y) := \sin^2 (2\pi y - \cos 2\pi x)
$$
generates a system for which $\Pfp \ne \P^f$.
This is because $f$ vanishes (even to second order) on a curve 
$$
y = {1 \over 2 \pi} \cos 2 \pi x + \Z
\label e.curve
$$
that is not in the zero set of any trigonometric polynomial in two
variables except the zero polynomial.
However, $f$ has the same distribution as the trigonometric polynomial
$g(x, y) := \sin^2 2 \pi y$ and
$(\P^g)^+ = \P^g$.
To see that the curve \ref e.curve/ does not lie in the zero set of any
nonzero trigonometric polynomial, suppose that $T(x, y)$ is a trigonometric
polynomial that vanishes there. Write $w := e^{2\pi i x}$ and $z := e^{2\pi
i y}$. By multiplying by a suitable complex
exponential, we may assume that $P(w, z) = T(x, y)$ is an analytic
polynomial. 
Our assumption is that 
$$
h(w) := P(w, e^{i(w+w^{-1})/2})
$$
satisfies $h(w) = 0$ for all $|w|=1$.
Since $h$ is an analytic function for $w \ne 0$, it follows that $h(w) = 0$
for all $w \ne 0$.
Now for each $z \ne 0$, there are an infinity of $w$ such that
$e^{i(w+w^{-1})/2} = z$. Therefore for each $z \ne 0$,
there are an infinity of $w$ such that $P(w, z) = 0$.
Since a polynomial has only a finite number of zeroes if it is not
identically zero, this means that for each $z \ne 0$, $P(w, z) = 0$ for all
$w \in \C$.
Hence $P \equiv 0$, as desired.
\endprocl

\procl x.algvar
For $d = 2$, if $f$ is real analytic on a neighborhood of its 
zero set, then $\Pfp = \P^f$ iff its zero set is contained in a 
(nontrivial) algebraic variety, where we view $\T^2$ as $ \{ (z_1, 
z_2) \in \C^2 \st |z_1| = |z_2| = 1 \}$.
This is because the slowest $f$ can vanish at a point is of order $x^2+y^2$
(since the constant and linear terms of $f$ must vanish)
and $1/(x^2+y^2)$ is not integrable. Therefore, all the zeroes of $f$ must be
cancelled by those of $T$.
Conversely, if the zero set of $f$ is contained in the zero set of a
trigonometric polynomial $T$, then by {\L}owasiewicz's inequality (see,
e.g., \ref b.BierMil/, Theorem 6.4), for a sufficiently large $n$, we have
$|T|^n/f$ is bounded, whence integrable.
\endprocl

\procl x.depends-distribution1
Here, we give a 1-dimensional example showing that the property $\Pfp =
\P^f$ does {\it not\/} depend only on the distribution of $f$, even among
continuous $f$.
Let $f(x) := \sqrt{x} |\sin (\pi/x)|$ on $[0, 1]$.
Then $\Pfp \ne \P^f$ since there are an infinite number of first-order zeroes.
However, let $g$ be the increasing rearrangement of $f$ on $[0, 1]$ and define
$$
h(x) := \cases{g(2x) &if $x \le 1/2$,\cr
               g(2-2x) &if $x \ge 1/2$.\cr}
$$
Then $h$ is continuous and has the same distribution as $f$, yet
$(\P^h)^+ = \P^h$ by \ref x.continuous-finite/ together with an easy
computation.
Of course, there is no such example for real-analytic $f$ in one dimension
by \ref x.continuous-finite/.
\endprocl

We finally state and prove our necessary and sufficient condition for
the strong full $K$ property.

\procl c.sn
Consider $f:\Td\to [0,1]$ and the corresponding probability measure
$\P^f$. Then
$\P^f$ is strong full $K$ if and only if there is a nonzero
trigonometric polynomial $T$ such that 
${|T|^2 \over f (\bfo-f)} \in L^1(\Td)$.
\endprocl

\proof
According to \ref t.sn/, 
$\Pfp=\P^f$ iff there exists a nonzero
trigonometric polynomial $T_1$ such that
${|T_1|^2 \over f} \in L^1(\Td)$. The same argument applied to $\bfo - f$
tells us that $\Pfm=\P^f$ iff there exists a nonzero
trigonometric polynomial $T_2$ such that
${|T_2|^2 \over \bfo-f} \in L^1(\Td)$. \ref l.pt/ now completes the proof.
(Take $T := T_1 T_2$.)
\Qed

\procl r.analytic
Observe that by trivial scaling, we could have an $f$ whose Fourier
coefficients go to zero slowly but which is bounded away from 0 and 1.
Hence slow decay of Fourier coefficients does not imply phase multiplicity.
On the other hand, if the coefficients decay exponentially, then 
$\sum_{n \in \Z} \widehat{f}(n) z^n$ is complex-analytic in an annulus around
the unit circle, which implies that its restriction to the unit circle is
real-analytic and hence there is phase uniqueness.
\endprocl

\procl r.folner
We mention that if $\Pfp=\P^f$, then this measure is
F{\o}lner independent in the sense of \ref b.Adams/ and
if in addition $\Pfm=\P^f$, then this measure is even strong
F{\o}lner independent in the sense that any conditioning outside
a box yields a measure which is $\dbar$ close to the unconditioned
process. The arguments for proving these facts are analogous to those
of \ref b.OW:ising/ (see \ref b.Adams/ for a published version), where it is
proved that the plus state for the Ising model is a Bernoulli shift.
(The concept of F{\o}lner independence has been used also in \ref b.HS/ and
\ref b.hoffman:MRF/.)
\endprocl  

We close this section with an interpretation of
our phase multiplicity results in terms 
of interpolation questions for wide-sense stationary 
processes. The following comments can be proved from \ref t.closest/ 
together with the well-known correspondence 
between interpolation and Szeg\H{o} infima.
First, $\Pfp= \delta_\bfz$ is equivalent to being able to 
interpolate perfectly from information far away for 
a wide-sense stationary process $\Seq{Y_n}$  with spectral density $f$,
which means that for every $n$, $Y_0$ is in the closed linear
subspace spanned by  $\{Y_k\}_{|k|\ge  n}$. For $d=1$, it
is stated on p.~102 of \ref b.Roz/ that this latter condition is 
equivalent to the negation of condition (iii) in \ref t.sn/.
A similar statement holds for $\Pfm$ with $f$ replaced by $\bfo-f$.
The proof of \ref t.sn/ implies (via this whole correspondence) that, 
for spectral measures that are
absolutely continuous with a bounded nonnegative density,
if perfect linear interpolation fails, then
our ability to interpolate as $n$ gets large goes to 0; i.e.,
the length of the projection of $Y_0$ onto the closed linear
subspace spanned by  $\{Y_k\}_{|k|\ge  n}$ goes to 0 as $n\to\infty$.

\bsection{1-Sided Phase Multiplicity}{s.1side}

We first study the notion of strong $K$ for $d=1$. 
This is natural since the
structure of ``one-sided'' behavior can be very different from
``two-sided'' behavior in various situations; for example, the existence
of bilaterally deterministic Bernoulli shift processes demonstrates this.
In contrast with \ref x.depends-distribution1/, we shall see that whether
$\P^f$ is strong $K$ depends only on the distribution of $f$.

Consider a function $f:\T\to [0,1]$ and the 
corresponding probability measure
$\P^f$ on $2^{\Z}$. We define a probability measure
$(\P^f)^{+,1}$ on $2^{\Z}$ thought of as ``$\P^f$ conditioned 
on all 1s at $-\infty$'' (the superscript ``1'' refers to the fact that we are
doing this on 1 side); namely,
$$
(\P^f)^{+,1}:= \lim_{n\to \infty} \lim_{k \to\infty} 
\P^f\big[\;\cdot \mid \eta\equiv 1 \hbox{ on }
[-n-k,-n]\cap \Z\big]
\,.
$$
The existence of the limit and its translation invariance follow from
stochastic monotonicity, as did that of $(\P^f)^+$.
As before, we also have $(\P^f)^{+,1} \preccurlyeq  \P^f$. 
In the analogous way (using 0 instead of 1 boundary conditions), one defines
$(\P^f)^{-,1}$, which then satisfies $\P^f  \preccurlyeq (\P^f)^{-,1}$. 

\procl d.pt1
The probability measure $\P^f$ has a {\bf 1-sided phase multiplicity} if 
$(\P^f)^{-,1}\neq(\P^f)^{+,1}$, and otherwise has {\bf 1-sided phase
uniqueness}.
\endprocl

Note the following analogue (whose proof is left to the reader) 
of \ref l.pt/ for the 1-sided case.

\procl l.pt1
$\P^f$ has 1-sided phase uniqueness if and only if it 
is strong $K$.
\endprocl

Our main result for 1-sided phase multiplicity is the following. 

\procl t.sn1GM
If $f:\T\to [0,1]$, then $(\P^f)^{+,1}=\P^f$ iff $\GM(f) > 0$. 
Moreover, if $(\P^f)^{+,1}\neq \P^f$, then $(\P^f)^{+,1}= \delta_\bfz$.
\endprocl

\proof
By \ref t.szego/, we have
$$
\lim_{k \to\infty} \P^f[\eta(0) = 1 \mid \eta_{-n} = \eta_{-n-1} =
\cdots = \eta_{-n-k} = 1]
=
\sum_{l=0}^{n-1} |\widehat{\varphi_f}(l)|^2
\,.
$$
(In fact, via \ref t.closest/, this special case of \ref t.szego/ is due to
Kolmogorov and Wiener; see \ref b.GrenanderSzego/, Section 10.9.)
Taking $n \to\infty$, we find that
$$
(\P^f)^{+,1}\big[\eta(0) = 1\big] = \|\widehat{\varphi_f}\|^2_2 = \|\varphi_f\|^2_2
\,.
$$
If $\GM(f) > 0$, then we obtain
$$
(\P^f)^{+,1}\big[\eta(0) = 1\big] = \|f\|_1 = \widehat
f(0) = \P^f\big[\eta(0) = 1\big]
\,.
$$
Since $(\P^f)^{+,1} \preccurlyeq \P^f$ and both probability measures are
$\Z$-invariant, it follows by \ref l.weak/ that $(\P^f)^{+,1} = \P^f$.
On the other hand, if $\GM(f) = 0$, then 
$$
(\P^f)^{+,1}\big[\eta(0) = 1\big] = 0
\,,
$$
so
$(\P^f)^{+,1} = \delta_\bfz$.
\Qed

Our necessary and sufficient condition for the strong $K$ property 
now follows immediately.

\procl c.sn1
For $f:\T\to [0,1]$, the corresponding probability measure $\P^f$
is strong $K$ if and only if $\GM(f) \GM(\bfo-f) > 0$.
\endprocl

We now construct a process $\P^f$ that is strong $K$ but not strong full $K$.

\procl t.diff
There exists an $f:\T\to [0,1]$ such that $\P^f$ is strong $K$ but not
strong full $K$.
\endprocl

\proof
Choose $f$ such that $f$ is continuous, bounded away from 1, vanishes at a
single point $x_0$, and such that $f(x) = 
e^{-1/(|x-x_0|)^{1/2}}$ in some neighborhood of $x_0$.
Since a trigonometric polynomial vanishes at its zeroes at most
polynomially quickly, there cannot exist a nonzero
trigonometric polynomial $T$
such that $|T|^2/f \in L^1(\T)$. Hence by \ref c.sn/,
$\P^f$ is not strong full $K$. 
On the other hand, it is clear that $\GM(f) \GM(\bfo-f) > 0$, 
so by \ref c.sn1/, we know that $\P^f$ is strong $K$.
\Qed

\procl r.plus
By combining Theorems 
\briefref t.dom/ and \briefref t.sn1GM/, we see that
$\P^f \succcurlyeq \P^g$ implies 
$(\P^f)^{+,1} \succcurlyeq (\P^g)^{+,1}$.
However, 
it is not necessarily the case that $(\P^f)^+ \succcurlyeq (\P^g)^+$.
For example, let $g \equiv 1/2$ and $f$ be a function as described in the last
result, but which has a geometric mean larger than $1/2$
(recall \ref t.dom/). 
This example also shows that if $\P^g \preccurlyeq \P^f$,
there do not necessarily exist $g'\le f'$ 
such that $\P^{g'} = \P^g$ and $\P^{f'} = \P^f$.  To see this,
let $f$ and $g$ be as above. Then $\P^{g'} = \P^g$ easily implies that
$g' \equiv 1/2$, which in turn (using $g'\le f'$) implies that
$(\P^{f'})^+ =  \P^{f'}$. Since
$\Pfp \neq  \P^f$, this contradicts the fact that $\P^f = \P^{f'}$.
\endprocl

It may be of interest to see another proof (which we only sketch)
of \ref t.sn1GM/ that does not
depend on the asymptotics of \ref t.szego/, but rather follows
the lines of the proof of \ref t.sn/.
The appropriate replacement of the trigonometric polynomials
for this question is the set 
$$
\calA :=
\{h\in L^2(\T) \st \hbox{ there exists } \ell  \hbox{ such that }
\widehat{h}(k) =0 \hbox{ for } k< \ell\}\,.
$$
Note that when $\ell$ is taken to be 0, we get the {\bf Hardy space}
of analytic functions, denoted $H^2(\T)$.

\procl t.sn1
Assume that $f:\T\to [0,1]$ and let $\calA$ be as defined as above.
The following are equivalent. 
\beginitems
\itemrm{(i)} $(\P^f)^{+,1}=\P^f$; 
\itemrm{(ii)} $f$ is in the closure in $L^2(1/f)$ of $\calA\cap L^2(1/f)$;
\itemrm{(iii)} There exists a nonzero element $T\in \calA$ such that 
${|T|^2 \over f}\in L^1(\T)$; i.e., $\calA\cap L^2(1/f)\neq {0}$;
\itemrm{(iv)} $\GM(f) > 0$. 
\enditems
Moreover, if $(\P^f)^{+,1}\neq \P^f$, then $(\P^f)^{+,1}= \delta_\bfz$.
\endprocl

\proof
{\bf (i) iff  (ii)} and {\bf (ii) implies (iii):} These are proved
exactly as in \ref t.sn/.

{\bf (iii) implies (ii):}
Fix a nonzero $T \in \calA$ with ${|T|^2 \over f}\in L^1(\T)$.
Given any $\epsilon$, there exists $\delta>0$ such that if $B$ is any
subset of $\T$ with measure less than $\delta$, then
$$
\int_B \left({|T|^2 \over f} + 2|T| +f\right) < \epsilon\,.
$$
Since $T$ vanishes only on a set of measure 0
(see \ref b.Rudin:RCA/, Theorem 17.18, p.~345),  
there is a $\gamma>0$ such that
$\{x\in \T \st |T(x)| < \gamma\}$ has measure less than $\delta$.
If we take $B$ to be this latter set, then we can proceed exactly 
as in \ref t.sn/ since now $T$ is bounded away from 0 on $B^c$. 

{\bf (iii) implies (iv):}
Any function in $\calA$ is a complex exponential times a function in the
Hardy space $H^2(\T)$. Hence there is 
a nonzero function $T \in H^2(\T)$ with $|T|^2/f \in
L^1(\T)$. Now
$\GM(|T|) > 0$  (see \ref b.Rudin:RCA/, Theorem 17.17, p.~344). 
Letting $g := |T|^2/f$, we have that
$g \in L^1(\T)$ and hence $\GM(g) \le \int g < \infty$.
Therefore $\GM(f) = \GM(|T|^2/g) = \GM(|T|)^2/\GM(g) > 0$, as desired.

{\bf (iv) implies (iii):}
If $\GM(f) > 0$, then take $T := \varphi_f$, defined in
\ref s.dom/.
For this $T$, we have $|T|^2/f = \bfo$, which is trivially integrable.

The final statement is proved as in \ref t.sn/.
\Qed

We now extend \ref t.sn1GM/ and \ref c.sn1/ to higher dimensions.
Fix a choice of ordering of $\Zd$ given by a set $\ord$ 
and a sequence $\Seq{k_n}\in \Z^d$ such that the distance from
$\bfz$ to $-(\ord+k_n)$ goes to $\infty$.

\procl d.1sK
A translation-invariant probability measure $\mu$ on $2^{\Z^d}$ 
{\bf has phase uniqueness 
relative to the ordering induced by $\ord$ and the sequence $\Seq{k_n}$ \/} 
as above if for any finite cylinder event $E$ and 
for all $\epsilon>0$ there exists an $N$ such that for all $n\ge N$
$$
\big|\mu\big(E \mid \F(-(\ord+k_n)\big)-\mu(E)\big| <\epsilon
\quad\mu\hbox{-a.e.}
$$
\endprocl

It is easy to check using \ref p.JNRD/
that this definition {\it does not depend\/} on the choice
of the sequence $\Seq{k_n}$. However,
it turns out that the choice of ordering, $\ord$, of $\Zd$ can make a 
difference. In particular, if the ordering is archimedean
(meaning that for any two positive elements $a$ and $b$, 
there exists an integer $n$ such that $na > b$), 
such as\ftnote{${}^*$}{In fact, all archimedean orders arise in this way; 
see \ref b.Teh/, \ref b.Zaiceva/, or \ref b.Trevisan/.} when the ordering is
induced by $\ord := \{ k \in \Z^d \st k \cdot x > 0 \}$, 
where $x \in \R^d$ is a fixed vector having two coordinates with an
irrational ratio, then we have a complete characterization
of 1-sided phase multiplicity in terms of the geometric mean.  For the standard
lexicographic ordering, however, 1-sided phase multiplicity cannot be
characterized by the geometric mean alone, as shown by the
example $f :=   \I{[0,1]\times[0, 1/2]}$ in two dimensions.
Here, both $f$ and $\bfo - f$ have 0 geometric mean, but the columns of a
configuration are independent under $\P^f$, so that conditioning on the remote
past has no effect on the present.
On the other hand, positivity of the geometric means will still {\it suffice\/}
for phase uniqueness 
relative to the ordering induced by $\ord$,
as we now prove. Analogously to the
one-dimensional case, we let, for our given ordering $\ord$,
$$
(\P^f)^{+,1}:= \lim_{n\to \infty}
\P^f\big[\;\cdot \mid \eta\equiv 1 \hbox{ on } -(\ord+k_n)\big]
\,.
$$
By \ref p.JNRD/, this limit exists and is independent of $\Seq{k_n}$.
As before, $(\P^f)^{-,1}$ is defined analogously using boundary
conditions of 0s, one-sided phase multiplicity is defined
by $(\P^f)^{-,1}\neq(\P^f)^{+,1}$, and the property
in \ref d.1sK/
is equivalent to one-sided phase uniqueness.

\procl t.sn1GM-d
Let $f:\Td\to [0,1]$. If $\GM(f) > 0$, then $(\P^f)^{+,1}=\P^f$.
If $\GM(f) = 0$ and the ordering is archimedean, then
$(\P^f)^{+,1}= \delta_\bfz$. In particular, if the 
ordering is archimedean, then $\P^f$ has phase 
uniqueness (relative to this ordering) if and only if 
$\GM(f) \GM(\bfo-f) > 0$.
\endprocl

We begin with the background results on which we rely for the proof.
For an ordering induced by $\ord$, define the 
corresponding {\bf Helson-Lowdenslager space}
$$
\hl^2 := \hl^2(\T^d, \ord) := \Big \{ f \in L^2(\T^d) \st \supp \widehat f
\subset \ord \cup \{ \bfz \} \Big \}
\,.
$$
For $0 \le f \in L^1(\lambda)$, let $[\ord \cup \{ \bfz \}]$ be the linear span
of $ \{ \ee_k \st k \in \ord \cup \{ \bfz \} \}$ and $[\ord \cup \{ \bfz
\}]_f$ be its closure in $L^2(f)$.
The replacement for outer functions is the class of {\bf spectral factors},
which are the functions $\varphi \in \hl^2$ with the additional properties 
$$
\widehat\varphi(\bfz) > 0
\label e.at0
$$
and
$$
1/\varphi \in [\ord \cup \{ \bfz \}]_{|\varphi|^2}
\,.
\label e.inverse
$$
\ref b.HelLow:I/ show that for any $\ord$ and for any
$0 \le f \in L^1(\T^d)$, the condition $\GM(f) >
0$ is equivalent to the existence of a spectral factor $\varphi_f$ such that
$|\varphi_f|^2 = f$ a.e. [More precisely, they prove $\GM(f) > 0$ iff $\exists
\varphi \in \hl^2$ satisfying $|\varphi|^2 = f$ a.e.\ and 
\ref e.at0/. Their proof
shows that in this case, $\varphi$ can be chosen so that also \ref e.inverse/
holds.] Furthermore, \ref b.Lyons:szego/ shows that $\varphi_f$ is then unique.

Denote by $\nu_{f}$ the measure on $2^{\ord \cup \{ \bfz \}}$ given
by $\P^f[\;\cdot \mid \eta \restrict (-\ord) \equiv 1]$. This is, as usual,
defined by conditioning on more and more 1s in $-\ord$ and using
\ref p.JNRD/. The results of \ref b.Lyons:szego/ imply the following.

\procl t.szego-d 
Fix an ordering induced by a set $\ord$.
Let $f : \Td \to [0, 1]$ be measurable. If $\GM(f) > 0$, define
$\varphi_f$ to be its spectral factor (with respect to $\ord$).
Otherwise, define $\varphi_f := \bfz$.
If $\GM(f) > 0$ or the ordering of $\Zd$ given by $\ord$ is archimedean, then
the measure $\nu_{f}$ is equal to
the determinantal probability measure corresponding to the positive
contraction on $\ell^2(\ord \cup \{ \bfz \})$ whose $(j, k)$-matrix entry is 
$$
\sum_{l\in \ord \cup \{ \bfz \},\; l \preccurlyeq j,k} 
\overline{\widehat{\varphi_f}(j-l)}\widehat{\varphi_f}(k-l)
\,.
$$
\endprocl

\proofof t.sn1GM-d
Suppose first that $\GM(f) > 0$ or the ordering given by $\ord$ is archimedean.
By \ref t.szego-d/, we have
$$
\P^f[\eta(\bfz) = 1 \mid \eta \restrict -(\ord+k_n) \equiv 1]
=
\P^f[\eta(k_n) = 1 \mid \eta \restrict (-\ord) \equiv 1]
=
\sum_{l \in \ord \cup \{ \bfz \},\; l \preccurlyeq k_n} 
|\widehat{\varphi_f}(k_n - l)|^2
\,.
$$
Taking $n \to\infty$, we find that
$$
(\P^f)^{+,1}\big[\eta(\bfz) = 1\big] = \|\widehat{\varphi_f}\|^2_2 =
\|\varphi_f\|^2_2
\,.
$$
If $\GM(f) > 0$, then we obtain
$$
(\P^f)^{+,1}\big[\eta(\bfz) = 1\big] = \|f\|_1 = \widehat
f(\bfz) = \P^f\big[\eta(\bfz) = 1\big]
\,.
$$
Since $(\P^f)^{+,1} \preccurlyeq \P^f$ and both probability measures are
$\Zd$-invariant, it follows that $(\P^f)^{+,1} = \P^f$.
However, if $\GM(f) = 0$ and the ordering is archimedean, then 
$$
(\P^f)^{+,1}\big[\eta(\bfz) = 1\big] = 0
\,,
$$
so
$(\P^f)^{+,1} = \delta_\bfz$.
The last statement of the theorem follows as before.
\Qed

Analogously to \ref s.sk/, our one-sided phase multiplicity results 
for $d=1$ can be translated to known results for prediction (rather than
interpolation) questions for wide-sense stationary 
processes. We leave these observations to the reader. 

\bsection{Open Questions}{s.open}

\procl q.1
Calculate $H(\P^f)$.
\endprocl

We conjecture that entropy is concave:

\procl g.concave
For any $f$ and $g$, we have $H\big(\P^{(f+g)/2}\big) \ge \big(H(\P^f) +
H(\P^g)\big)/2$.
\endprocl

\procl q.6
Letting $A$ vary over all sets of measure 1/2, how do we get the largest
and the smallest entropies for $\P^{\I A}$?
More generally, given $f$, which $g$ with the same distribution as $f$
maximize or minimize $H(\P^g)$?
\endprocl

\procl q.bounds 
If $d=1$ and $\GM(f) \GM(\bfo - f) > 0$, do we get arbitrarily close lower
bounds (in principle) on $H(\P^f)$ by the method in \ref s.dom/?
How can one get close lower bounds on the entropy for higher-dimensional
processes?
\endprocl

\procl q.block
Suppose $f : \T \to [0, 1]$ is a trigonometric polynomial $f$ of
degree $m$.
Then $\P^f$ is $m$-dependent, as are all $(m+1)$-block factors of
independent processes (as defined by \ref b.Hed/).
Is it the case that $\P^{f}$ is an $(m+1)$-block factor of an i.i.d.\ process?
Erik Broman (personal communication) has shown this when $m = 1$.
If one can find such block factors sufficiently explicitly for
trigonometric polynomials, then one could find explicit factors of i.i.d.\
processes that give any process $\P^f$. This would enable one to use more
standard probabilistic techniques to study $\P^f$.
More generally, in higher dimensions,
if $f : \Td \to [0, 1]$ is a trigonometric polynomial, is $\P^f$ a block
factor?
\endprocl

\procl q.otherPT
Given an ordering $\ord$ on $\Zd$ and an increasing sequence of finite sets
$\ord_n \subset \ord$ whose union is all of $\ord$, when is there phase
multiplicity in the sense that 
$$
\lim_{n\to\infty}
\P^f[\eta(\bfz) = 1 \mid \eta \restrict (\ord \setminus \ord_n) \equiv 1]
\ne
\lim_{n\to\infty}
\P^f[\eta(\bfz) = 1 \mid \eta \restrict (\ord \setminus \ord_n) \equiv 0]
\,?
$$
This clearly does not depend on the choice of $\ord_n$.
In one dimension, this is the same as 1-sided phase multiplicity, where we know
the answer.
\endprocl

\procl q.3
Suppose that $d = 1$.
Note that translation and flip of $f$ yield the same measure $\P^{f}$.
Does $\P^{f}$ determine $f$ up to translation and
flip? 
\endprocl

We now ask about some properties that are important for models
in statistical physics.
Let $\others$ be as in \ref d.fulldom/.
A stationary process $\nu$ is called {\bf quasi-local} if $\nu[\eta(\bfz) =
1 \mid \others]$ has a continuous version, where the product topology is
used on $2^{\Zd \setminus \{ \bfz \}}$.
It is called {\bf almost surely quasi-local} if $\nu[\eta(\bfz) =
1 \mid \others]$ has an almost surely continuous version.

\procl q.q-local
For which $f$ is $\P^f$ quasi-local or almost surely quasi-local?
When there is phase multiplicity, then
$\P^f$ is as far as possible from this since then each version of
$\P^f[\eta(\bfz) = 1 \mid \others]$ is nowhere continuous.
In one dimension,
there is a natural definition of one-sided quasi-locality, which is the same
as what is sometimes referred to as the uniform martingale property
(see \ref b.kalikow/). It is
not hard to see that this property, or even its almost sure version, would
imply an affirmative answer to \ref q.bounds/.
\endprocl

We say that a stationary process $\nu$ on $2^\Zd$ is {\bf insertion
tolerant} if $\nu[\eta(\bfz) = 1 \mid \others] > 0$ $\nu$-a.s.
Similarly, $\nu$ is {\bf deletion tolerant} if $\nu[\eta(\bfz) = 0 \mid
\others] > 0$ $\nu$-a.s.

\procl g.tolerance
Every $\P^f$ is insertion and deletion tolerant, other than the trivial cases
$f = \bfz$ and $f = \bfo$.
\endprocl

\ref t.fulldom/ implies a positive solution to \ref g.tolerance/ when
$1/[f(\bfo-f)]$ is integrable. Indeed, we see by \ref t.fulldom/ that this is
the exact criterion for uniform insertion tolerance and uniform
deletion tolerance.

\medbreak
\noindent {\bf Acknowledgements.}\enspace
We thank Chris Hoffman for his help with \ref p.hoffman/.
We thank Eric Bedford and Jeff Geronimo for important references and Jean-Paul
Thouvenot for discussions.
We are grateful to David Wilson for permission to include Figures \briefref
f.window-ust/ and \briefref
f.window-tiling/ and to Peter Duren and J.M.~Landsberg for helpful remarks.
\ref f.window-tiling/ was created using the linear algebraic techniques of
\ref b.wilson:planar/ for generating domino tilings; the needed
matrix inversion was
accomplished using the formulas of \ref b.Kenyon:local/. The resulting tiling
gives dual spanning trees by the bijection of Temperley (see \ref
b.KPW:tr-mat/). One of the trees is \ref f.window-ust/.

\bibfile{\jobname}
\def\noop#1{\relax}
\input \jobname.bbl

\filbreak
\begingroup
\eightpoint\sc
\parindent=0pt\baselineskip=10pt

Department of Mathematics,
Indiana University,
Bloomington, IN 47405-5701
\emailwww{rdlyons@indiana.edu}
{http://php.indiana.edu/\string~rdlyons/}

and

School of Mathematics,
Georgia Institute of Technology,
Atlanta, GA 30332-0160
\email{rdlyons@math.gatech.edu}

Department of Mathematics,
Chalmers University of Technology,
S-41296 Gothenburg, Sweden 
\email{steif@math.chalmers.se}

and

School of Mathematics,
Georgia Institute of Technology,
Atlanta, GA 30332-0160
\email{steif@math.gatech.edu}

\endgroup

\bye